\documentclass[11pt]{article}
\usepackage{amsmath}
\usepackage{amssymb, amsmath}
\let\pa=\partial
\let\al=\alpha

\let\f=\frac
\let\p=\psi

\def\no{\noindent}
\def\na{\nabla}

\def\p{\partial}

\def\dv{\mbox{div}}

\def\C{\mathop{\bf C\kern 0pt}\nolimits}
\def\DD{\mathop{\bf D\kern 0pt}\nolimits}
\def\K{\mathop{\bf K\kern 0pt}\nolimits}
\def\N{\mathop{\bf N\kern 0pt}\nolimits}
\def\Q{\mathop{\bf Q\kern 0pt}\nolimits}
\def\R{\mathop{\bf R\kern 0pt}\nolimits}

\newcommand{\ds}{\displaystyle}

\newcommand{\Z}{{\mathbf Z}}

\newcommand{\ef}{ \hfill $ \blacksquare $ \vskip 3mm}

\newcommand{\beq}{\begin{equation}}
\newcommand{\eeq}{\end{equation}}
\newcommand{\ben}{\begin{eqnarray}}
\newcommand{\een}{\end{eqnarray}}
\newcommand{\beno}{\begin{eqnarray*}}
\newcommand{\eeno}{\end{eqnarray*}}

\newtheorem{Def}{Definition}[section]

\newtheorem{rmk}{Remark}[section]

\renewcommand{\theequation}{\thesection.\arabic{equation}}

\newtheorem{Theorem}{Theorem}[section]

\newtheorem{Proposition}[Theorem]{Proposition}
\newtheorem{Lemma}[Theorem]{Lemma}

\begin{document}
\title{On the well-posedness for the viscous shallow water equations}
\author{Qionglei Chen $^\dag$  Changxing Miao $^{\dag}$ and Zhifei Zhang
$^{\ddag}$\\[2mm]
{\small $ ^\dag$ Institute of Applied Physics and Computational Mathematics, Beijing 100088, China}\\
{\small E-mail: chen\_qionglei@iapcm.ac.cn  and   miao\_changxing@iapcm.ac.cn}\\[2mm]
{\small $ ^\ddag$ School of  Mathematical Science, Peking University, Beijing 100871, China}\\
{\small E-mail: zfzhang@math.pku.edu.cn}}

\date{}
\maketitle

\begin{abstract}
In this paper, we prove the existence and uniqueness of the solutions for the 2D viscous shallow water equations
with low regularity assumptions on the initial data as well as the initial height bounded away from zero.
\end{abstract}

\renewcommand{\theequation}{\thesection.\arabic{equation}}
\setcounter{equation}{0}
%%%%%%%%%%%%%%%%%%%%%%%%%%%%%%%%%%%%%%%%%%%%%%
%%%%%%%%%%%%%%%%%%%%%%%%%%%%%%%%%%%%%%%%%%
\section{Introduction}

In this paper, we study the 2D viscous
shallow water equations with a more general diffusion
\begin{equation}\label{1.1}
\left\{
\begin{array}{ll}
h_t+\dv(hu)=0,\\
h(u_t+u\cdot\nabla u)-\nu\nabla\cdot(hD(u))-\nu\nabla(h \textrm{div}(u))+h\nabla h=0, \\
u(0,\cdot)=u_0,\,\,h(0,\cdot)=h_0 ,
\end{array}
\right.
\end{equation}
where $h(t,x)$ is the height of fluid surface,
$u(t,x)=(u_1(t,x),u_2(t,x))$ is the horizontal velocity vector
field, $D(u)=\f12(\na u+\na u^t)$ is the deformation tensor, and
$\nu>0$ is the viscous coefficient. If the diffusion terms in (\ref{1.1}) are
replaced by $-\nu\nabla\cdot(h\na u)$,
then (\ref{1.1}) turns into the usual viscous shallow water equations.

Recently, the viscous shallow water equations have been widely
studied by Mathematicians, see the review paper\cite{BDM}.  Bui
\cite{Bui} proved the local existence and uniqueness of classical
solutions to the Cauchy-Dirichlet problem for the shallow water
equations with initial data $h_0$, $u_0$ in H\"{o}lder spaces as
well as $h_0$ bounded away from vacuum. Kloeden\cite{Klo} and
Sundbye\cite{Su1} independently proved global existence and
uniqueness of classical solutions to the Cauchy-Dirichlet problem in
Sobolev spaces. Later, Sundbye\cite{Su2} also proved global
existence and uniqueness of classical solutions  to the Cauchy
problem. However, for all above results (except \cite{Bui}), the
authors only consider the case when the initial data $h_0$ is a
small perturbation of some positive constant $\bar h_0$ and $u_0$ is
small in some sense. Very recently, Wang and Xu\cite{WX} proved the local
well-posedness of the Cauchy problem in Sobolev spaces
for the large data $u_0$ and $h_0$ closing to $\bar h_0$. More
precisely, they obtained the following result.
\begin{Theorem}\label{WX}\cite{WX} Let $\bar{h}_0$ be a strictly positive constant and $s>2$.
Assume that
\beno&&(i)\quad(u_0,h_0-\bar{h}_0)\in
H^s(\R^2)\otimes H^s(\R^2);\hspace{7cm}\\
&&(ii)\quad\|h_0-\overline{h}_0\|_{H^s}\ll\overline{h}_0.\eeno Then
there exist a positive time $T$ and a unique solution $(u,h)$ of
$(\ref{1.1})$ such that
\ben\label{1.2} u,\,\, h-\bar h_0\in
L^\infty([0,T],H^s),\quad \nabla u\in L^2([0,T]; H^s).
\een
Moreover, there exists a strictly positive constant $c$ such that if
\begin{eqnarray}\label{1.3}
\|u_0\|_{H^s}+ \|h_0-\bar{h}_0\|_{H^s}\le c,
\end{eqnarray}
then we can choose $T=+\infty$.
\end{Theorem}

One purpose of this paper is to study the well-posedness of
$(\ref{1.1})$ for the initial data with the minimal regularity.
For the incompressible Navier-Stokes equations, such research has been initiated by
Fujita and Kato\cite{FK}, see also \cite{Can1, Can2, Mey2} for other relevant results. They proved local well-posedness for the incompressible Navier-Stokes equations
in the scaling invariant space. The scaling invariance means that
if $(u,p)$ is a solution of the incompressible Navier-Stokes equations with  initial data $u_0(x)$, then
\ben
\label{1.4}u_\lambda(t,x)\triangleq\lambda u(\lambda^2 t,\lambda x),\quad
p_\lambda(t,x)\triangleq \lambda^2 p(\lambda^2 t,\lambda x)
\een
 is also a solution of the incompressible Navier-Stokes equations with
$u_{0,\lambda}\triangleq \lambda u_0(\lambda x)$.
Obviously, $\dot H^{\frac d 2-1}(\R^d)$
is a scaling invariant space under the scaling of (\ref{1.4}), i.e.
$$\|u_\lambda\|_{\dot{H}^{\frac d 2-1}}=\|u\|_{\dot{H}^{\frac d 2-1}}.$$
The equations $(\ref{1.1})$ have no scaling invariance like the incompressible Navier-Stokes equations.
However, due to the similarity of the structure between (\ref{1.1}) and the incompressible Navier-Stokes equations,
we still solve (\ref{1.1}) for initial data whose regularity fits with the scaling of (\ref{1.4}).
It should be pointed out that R. Danchin was the first to consider the similar problem for the compressible
Navier-Stokes equations, and some ideas of this paper is motivated by \cite{D1}.

The second purpose of this paper is to prove the local well-posedness of (\ref{1.1}) under  more natural assumption that
the initial height is bounded away from zero.
For the initial data with slightly higher regularity,
this can be easily obtained by modifying the argument of Danchin\cite{D5}.
However, for the initial data with low regularity, his method is not applicable any more, since
the proof of \cite{D5} relies on the fact that some profits can be gained from the inclusion map
$B^s\hookrightarrow L^\infty$ in the case of $s>\frac d2$.
For this reason, we have to introduce some kind of weighted Besov space $E^s_T$(see Section 3) which is crucial to
get rid of the condition that the initial height $h_0$ is close to $\bar{h}_0$.
One important observation is that the $E^s_T$ norm of the solution is small for small time $T$.

Before stating our main result, let us first introduce some notations and definitions.
Choose a radial function $\varphi\in{\cal S}(\R^d)$ such that

$${\rm supp}\,\varphi\subset\{\xi\in \R^d;
\frac{5}{6}\le|\xi|\le \frac{12}{5}\},\quad \sum_{k\in\Z}\varphi(2^{-k}\xi)=1,
\quad\xi\in\R^d\setminus \{0\}.$$
Here $\varphi_k(\xi)=\varphi(2^{-k}\xi)$,
$k\in \Z$.

\begin{Def}\label{Def1.1}
Let $k\in \Z$, the Littlewood-Paley projection operators $\Delta_k$ and $S_k$ are defined as follows
$$\Delta_kf=\varphi(2^{-k}D)f,\quad
S_kf=\sum_{j\le k-1}\Delta_jf, \quad \textrm{for}\quad f\in {\cal S'}(\R^d).$$
\end{Def}

We denote the space ${\cal Z'}(\R^d)$ by the dual space of ${\cal
Z }(\R^d)=\{f\in {\cal S}(\R^d);\,D^\alpha \hat{f}(0)=0;
\forall\alpha\in\N^d \,\mbox{multi-index}\}$, it also can be
identified by the quotient space of ${\cal S'}(\R^d)/{\cal P}$ with the
polynomials space ${\cal P}$. The formal equality $$f=\sum_{k\in\Z}\Delta_kf$$ holds true for $f\in {\cal Z'}(\R^d)$
and is called the homogeneous Littlewood-Paley decomposition. It
has nice properties of quasi-orthogonality: with our choice of
$\varphi$,
\begin{eqnarray}\label{1.5}
\Delta_j\Delta_kf=0\quad i\!f\quad|j-k|\ge 2\quad \textrm{and}
\quad \Delta_j(S_{k-1}\Delta_k
f)=0\quad i\!f\quad|j-k|\ge 4.
\end{eqnarray}

\begin{Def}\label{Def1.2}
Let $s\in\R$, $1\le p, r\le+\infty$. The homogeneous Besov space $\dot{B}^{s}_{p,r}$
is  defined by
$$\dot{B}^{s}_{p,r}=\{f\in {\cal Z'}(\R^d):\,\|f\|_{\dot{B}^{s}_{p,r}}<+\infty\},$$
where
\beno
&&\|f\|_{\dot{B}^{s}_{p,r}}=\left\{\begin{array}{l}\displaystyle\bigg(\sum_{k\in\Z} 2^{ksr}
\|\Delta_kf\|^r_{p}\bigg)^{\frac{1}{r}},\quad\mbox{for}\quad r<+\infty,\\
\displaystyle\sup_{k\in\Z}2^{ks}\|\Delta_kf\|_{p},\quad\mbox{for}\quad r=+\infty.\end{array}\right.\eeno
\end{Def}
If\, $p=r=2$, $\dot{B}^s_{2,2}=\dot{H}^s$,
and if\, $d=2$, we have $\dot{B}^{1}_{2,1}\hookrightarrow L^\infty$ and
$$\|f\|_{\infty}\le C\|f\|_{\dot{B}^{1}_{2,1}}.$$
We refer to \cite{Ch1,Tr} for more details.

In addition to the general time-space space such as $L^\rho(0,T; \dot{B}^{s}_{p,r})$,
we introduce a useful mixed time-space homogeneous Besov space
$\widetilde{L}^\rho_T(\dot{B}^{s}_{p,r})$
which is initiated in \cite{CL} and is used in the proof of the uniqueness.
\begin{Def}\label{Def1.3}
Let $s\in\R$, $1\le p,r,\rho\le+\infty$, $0<T\le+\infty$. The  mixed time-space homogeneous
Besov space $\widetilde{L}^r_T(\dot{B}^s_{p,q})$ is defined by
$$\widetilde{L}^\rho_T(\dot{B}^{s}_{p,r})=
\{f\in {\cal Z'}(\R^{d+1}):\,\|f\|_{\widetilde{L}^\rho_T(\dot{B}^{s}_{p,r})}<+\infty\},$$\
where$$\|f\|_{\widetilde{L}^\rho_T(\dot{B}^{s}_{p,r})}=
\bigg\|2^{ks}
\bigg(\int_0^T\|\Delta_kf(t)\|^\rho_{p}dt\bigg)
^{\frac{1}{\rho}}\bigg\|_{\ell^r}.$$
\end{Def}
Using the Minkowski inequality, it is easy to verify that
$$L^\rho_T(\dot{B}^{s}_{p,r})\subseteq\widetilde{L}^\rho_T(\dot{B}^{s}_{p,r})
\quad\hbox{if}\quad \rho\le r\quad \textrm{and}\quad
\widetilde{L}^\rho_T(\dot{B}^{s}_{p,r})\subseteq{L}^\rho_T(\dot{B}^{s}_{p,r})
\quad\hbox{if}\quad \rho\ge r.$$

Next, we introduce a hybrid-index Besov space
which plays an important role in the study of compressible
fluids and is initiated in \cite{D1,D2}.
\begin{Def}\label{Def1.4}
Let $s$, $\sigma\in\R$, and set
$$\|f\|_{\widetilde{B}^{s,\sigma}_2}\triangleq\sum_{k\le0}2^{ks}\|\Delta_k f\|_{2}
+\sum_{k>0}2^{k\sigma}\|\Delta_k f\|_{2}.$$
Let $m=-[\frac{d}{2}+1-s]$, we define
\beno
&&\widetilde{B}^{s,\sigma}_2(\R^d)=\{f\in{\cal S'}(\R^d):
\|f\|_{\widetilde{B}^{s,\sigma}_2}<+\infty\}\quad
\hbox{if}\quad m<0,\\
&&\widetilde{B}^{s,\sigma}_2(\R^d)=\{f\in{\cal S'}(\R^d)/{\cal P}_m:
\|f\|_{\widetilde{B}^{s,\sigma}_2}<+\infty\}\quad
\hbox{if}\quad m\ge0,\eeno
where ${\cal P}_m$ denotes the set of polynomials of degree $\le m$.
\end{Def}

Throughout this paper, we will denote
$\dot{B}^{s}_{2,1}$ by $B^s$, and $\widetilde{B}^{s,\sigma}_2$ by
$\widetilde{B}^{s,\sigma}$. The following facts can be easily verified by using the definition of $\widetilde{B}^{s,\sigma}$:\vspace{0.1cm}

(i) $\widetilde{B}^{s,s}=\dot{B}^s_{2,1}$;

(ii)\,If $s\le \sigma$, then $\widetilde{B}^{s,\sigma}=\dot{B}^s_{2,1}\cap\dot{B}^\sigma_{2,1}$.
Otherwise, $\widetilde{B}^{s,\sigma}=\dot{B}^s_{2,1}+\dot{B}^\sigma_{2,1}$.
\vspace{0.2cm}

Now we state our main result as follows.
\begin{Theorem}\label{Them1.2} Let $\bar{h}_0$ be a positive constant.
Assume that

(i)\, $(u_0,h_0-\bar{h}_0)\in B^0(\R^2)\otimes
\widetilde{B}^{0,1}(\R^2);$\vspace{0.1cm}

(ii)\,\,$h_0\ge \bar{h}_0$.\\ \vspace{0.1cm}
Then there exist a positive time $T$ and a unique   solution $(u,h)$ of $(\ref{1.1})$ such that
\ben\label{1.6}
u\in C([0,T]; B^0)
\cap L^1(0,T; B^2),\,\,
h-\bar{h}_0\in C([0,T]; \widetilde{B}^{0,1})\cap
L^1(0,T; \widetilde{B}^{2,1}),\quad h\ge \f12 \bar{h}_0.
\een
Moreover, there exists a strictly positive constant $c$ such that if
\begin{eqnarray}\label{1.7}
\|u_0\|_{B^0}+
\|h_0-\bar{h}_0\|_{\widetilde{B}^{0,1}}\le c,\end{eqnarray}
then we can choose $T=+\infty$.
\end{Theorem}

The structure of this paper is as follows.

In Section 2, we recall some useful multilinear estimates in the Besov spaces .
In Section 3, we prove the existence of solution. In Section 4,
we prove the uniqueness of the solution.
Finally, in the Appendix, we prove some multilinear estimates in the weighted Besov spaces.

Throughout the paper, $C$ denotes various ``harmless'' large
finite constants, and $c$ denotes various ``harmless" small
constants. We shall sometimes use $X\lesssim Y$ to
denote the estimate $X\le CY$ for some constant $C$.
We denote $\|\cdot\|_p$ by the $L^p$ norm of a function.
\vspace{.3cm}

\setcounter{equation}{0}
\section{Multilinear estimates in the Besov spaces}

Let us first recall the Bony's paraproduct decomposition.

\begin{Def}\label{Def2.1}
We shall use the following  Bony's paraproduct decomposition(see \cite{BC,B})
\begin{eqnarray}\label{2.1}
fg=T_fg+T_gf+R(f,g),
\end{eqnarray}
with
\begin{eqnarray}\label{2.2}
T_fg=\sum_{k\in\Z}S_{k-1}f\Delta_kg\quad \textrm{and} \quad
R(f,g)=\sum_{k\in\Z}\sum_{|k'-k|\le 1}\Delta_{k} f\Delta_{k'}g.
\end{eqnarray}
\end{Def}

Next, let us recall some useful lemmas and multilinear estimates in the Besov spaces.

\begin{Lemma}(Bernstein's inequality)\label{Lem2.1}
Let $1\le p\le q\le+\infty$. Assume that $f\in {\cal S'}(\R^d)$,
then for any $\gamma\in\Z^d$,
there exist
constants $C_1$, $C_2$ independent of $f$, $j$ such that
\beno
&&{\rm supp}\hat f\subseteq \{|\xi|\le
A_02^{j}\}\Rightarrow
\|\partial^\gamma f\|_q\le C_12^{j{|\gamma|}+j d(\frac{1}{p}-\frac{1}{q})}\|f\|_{p},
\\
&&{\rm supp}\hat f\subseteq \{A_12^{j}\le|\xi|\le A_22^{j}\}\Rightarrow
\|f\|_{p}\le C_22^{-j|\gamma|}\sup_{|\beta|=|\gamma|}\|\partial^\beta f\|_p.
\eeno
\end{Lemma}
The proof can be found in \cite{Ch1}.

\begin{Proposition}\label{Prop2.2}
If $s>0$, $f,g\in B^s\cap L^\infty$.
Then $fg\in B^s\cap L^\infty$ and
\beq\label{2.3}
\|fg\|_{B^s}\le C(\|f\|_{\infty}\|g\|_{B^s}+\|g\|_{\infty}\|f\|_{B^s}).
\eeq
If $s_1$, $s_2\le\frac{d}{2}$ such that $s_1+s_2>0$, $f\in B^{s_1}$, and
$g\in B^{s_2}$. Then $fg\in B^{s_1+s_2-\frac{d}{2}}$ and
\ben\label{2.4}
\|fg\|_{B^{s_1+s_2-\frac{d}{2}}}\le C\|f\|_{B^{s_1}}\|g\|_{B^{s_2}}.
\een
If $|s|<\frac{d}{2}$, $1\le r\le+\infty$, $f\in \dot{B}^s_{2,r}$ ,and
$g\in {B}^{\frac{d}{2}}$. Then $fg\in \dot{B}^{s}_{2,r}$ and
\beq\label{2.5}
\|fg\|_{\dot{B}^{s}_{2,r}}\le C\|f\|_{\dot{B}^{s}_{2,r}}
\|g\|_{{B}^{\frac{d}{2}}}.
\eeq
If $s\in(-\frac{d}{2},\frac{d}{2}]$, $f\in B^s$, and
$g\in \dot{B}^{-s}_{2,\infty}$. Then $fg\in \dot{B}^{-\frac{d}{2}}_{2,\infty}$ and
\ben\label{2.6}
\|fg\|_{\dot{B}^{-\frac{d}{2}}_{2,\infty}}\le C\|f\|_{{B}^{s}}
\|g\|_{\dot{B}^{-s}_{2,\infty}}.\een
If $1\le \rho_1, \rho_2, \rho\le \infty, s\in(-\frac{d}{2},\frac{d}{2}]$,
$f\in \widetilde{L}^{\rho_1}_T(B^s)$, and
$g\in \widetilde{L}^{\rho_2}_T(\dot{B}^{-s}_{2,\infty})$. Then there holds
\ben\label{2.7}
\|fg\|_{\widetilde{L}^\rho_T(\dot{B}^{-\frac{d}{2}}_{2,\infty})}\le
C\|f\|_{\widetilde{L}^{\rho_1}_T({B}^{s})}
\|g\|_{\widetilde{L}^{\rho_2}_T(\dot{B}^{-s}_{2,\infty})},
\een
where $\frac{1}{\rho_1}+\frac{1}{\rho_2}=\frac{1}{\rho}$.
\end{Proposition}
{\it Proof.}\quad
For the sake of simplicity, we only present the proof of (\ref{2.4}) below, the others can be
deduced in the same way (see also \cite{D3,RS}).
By the  Bony's paraproduct decomposition  and the
property of quasi-orthogonality (\ref{1.5}), for fixed $j\in\Z$, we write
\begin{align}
\Delta_j(fg)&=\sum_{|k-j|\le3}\Delta_j(S_{k-1}f\Delta_kg)+
\sum_{|k-j|\le3}\Delta_j(S_{k-1}g\Delta_kf)
+\sum_{k\ge j-2}\sum_{|k-k'|\le 1}\Delta_j(\Delta_{k}f\Delta_{k'}g)\nonumber\\
&\triangleq {I+II+III}.\nonumber\end{align}
Thanks to the definition of Besov space $B^s$, we have
\begin{align}\label{2.8}
\|fg\|_{B^{s_1+s_2-\frac{d}{2}}}&\le\bigg(\sum_{j\in
\Z}2^{(s_1+s_2-\frac{d}{2})j}\|I\|_{2}\bigg)
+\cdots+\bigg(\sum_{j\in
\Z}2^{(s_1+s_2-\frac{d}{2})j}\|III\|_{2}\bigg)\nonumber\\
&\triangleq {I'+II'+III'}.\end{align}
It suffices to  estimate the
above three terms separately.
Using the Young's inequality and lemma \ref{Lem2.1}, we have
\begin{align}
\|\Delta_j(S_{k-1}f\Delta_kg)\|_{2}&\lesssim
\|S_{k-1}f\|_{\infty}\|\Delta_kg\|_{2}
\lesssim
\sum_{k'\le k-2}\|\Delta_{k'} f\|_{\infty}\|\Delta_kg\|_{2}\nonumber\\
&\lesssim \sum_{k'\le k-2}2^{k' s_1}
\|\Delta_{k'} f\|_{2}2^{k'(\frac{d}{2}-s_1)}\|\Delta_kg\|_{2}
\nonumber\\
&\lesssim \|f\|_{B^{s_1}}
\|\Delta_kg\|_{2}2^{k(\frac{d}{2}-s_1)},\nonumber
\end{align}
where  we have used the fact $s_1\le \frac{d}{2}$ in the last inequality.
Hence, we get
\begin{align}\label{2.9}
{I'}&\lesssim \|f\|_{B^{s_1}}\sum_{j\in\Z}2^{(s_1+s_2-\frac{d}{2})j}
\sum_{|k-j|\le3}2^{k(\frac{d}{2}-s_1)}
\|\Delta_kg\|_{2}
\nonumber\\&\lesssim \|f\|_{B^{s_1}}
\sum_{|\ell|\le3}2^{-(s_1+s_2-\frac{d}{2})\ell}
\sum_{j\in\Z}2^{s_2(j+\ell)}
\|\Delta_{j+\ell}g\|_{2}
\lesssim \|f\|_{B^{s_1}}\|g\|_{B^{s_2}}.
\end{align}
Similarly, using the fact $s_2\le\frac{d}{2} $, we can obtain
\begin{align}\label{2.10}
{II'}&\lesssim \|f\|_{B^{s_1}}\|g\|_{B^{s_2}}.
\end{align}
Now we turn to estimate $III'$. From Lemma \ref{Lem2.1} and H\"{o}lder inequality
, it follows that
\beno
\|\Delta_j(\Delta_{k}f\Delta_{k'}g)\|_{2}\lesssim
2^{j\frac{d}{2}}\|\Delta_{k}f\Delta_{k'}g\|_{1}
\lesssim 2^{j\frac{d}{2}}\|\Delta_{k}f\|_{2}\|\Delta_{k'}g\|_{2}.
\eeno
So, we get by Minkowski inequality that for $s_1+s_2>0$
\begin{align}\label{2.11}
{III'}&\lesssim \sum_{j\in\Z}2^{(s_1+s_2-\frac{d}{2})j}2^{j\frac{d}{2}}
\bigg(\sum_{k\ge j-2}\sum_{|k-k'|\le1}
\|\Delta_{k}f\|_{2}
\|\Delta_{k'}g\|_{2}\bigg)
\nonumber\\&\lesssim
\sum_{\ell\ge-2}2^{-(s_1+s_2)\ell}
\sum_{j\in\Z}2^{s_1(j+\ell)}
\|\Delta_{j+\ell}f\|_{2}\|g\|_{B^{s_2}}
\lesssim \|f\|_{B^{s_1}}\|g\|_{B^{s_2}}.
\end{align}
Summing up (\ref{2.8})-(\ref{2.11}), we get the desired inequality
(\ref{2.4}). \ef

\begin{Proposition}\label{Prop2.3}
\textrm{(1)} Let $s>0$. Assume that $F\in
W^{[s]+2,\infty}_{loc}(\R^d)$ such that $F(0)=0$. Then there
exists a constant $C(s,d,F)$ such that if $u\in B^s\cap L^\infty$,
there holds \ben\label{2.12} &&\|F(u)\|_{B^s}\le
C(1+\|u\|_{\infty})^{[s]+1}\|u\|_{B^s} ;\een and if $u\in
\dot{B}^s_{2,\infty}\cap L^\infty$, there holds \ben\label{2.13}
&&\|F(u)\|_{\dot{B}^s_{2,\infty}}\le
C(1+\|u\|_{\infty})^{[s]+1}\|u\|_{\dot{B}^s_{2,\infty}} . \een
(2)Assume that  $G\in W^{[\frac{d}{2}]+3,\infty}_{loc}(\R^d)$ such
that $G'(0)=0$. Then there exists a functions $C(s,d,G)$ such that
if  $-\frac{d}{2}<s\le\frac{d}{2}$, $u, v\in B^{\frac{d}{2}}\cap
L^\infty$ and $u-v\in B^s$, there holds
\begin{align}\label{2.14}
\|G(u)-G(v)\|_{B^s}\le C(\|u\|_{\infty},\|v\|_{\infty})
(\|u\|_{B^{\frac{d}{2}}}+\|v\|_{B^{\frac{d}{2}}})
\|u-v\|_{B^s};\end{align}
and if $|s|<\frac{d}{2}$, $u, v\in B^\frac{d}{2}\cap L^\infty$ and $u-v\in\dot{B}^s_{2,\infty}$, there holds
\ben\label{2.15}
\|G(u)-G(v)\|_{\dot{B}^s_{2,\infty}}\le C(\|u\|_{\infty},\|v\|_{\infty})
(\|u\|_{B^{\frac{d}{2}}
}+\|v\|_{B^{\frac{d}{2}}})
\|u-v\|_{\dot{B}^s_{2,\infty}}
.\quad
\een
\end{Proposition}
{\it Proof.}\quad  We can refer to \cite{BC1,RS} for the proof of (1).
For (2), we refer to \cite{D1,D4}. For example, we write
$$G(u)-G(v)=(u-v)\int_0^1G'(v+\tau(u-v))d\tau,$$
then it follows from (\ref{2.5}) that for $|s|<\frac{d}{2}$
\beno
\|G(u)-G(v)\|_{\dot{B}^s_{2,\infty}}\le C\|u-v\|_{\dot{B}^s_{2,\infty}}
\|G'(v+\tau(u-v))\|_{B^\frac{d}{2}},
\eeno
which together with (\ref{2.12}) implies (\ref{2.15}).\ef

\begin{Proposition}\label{Prop2.4}
Let $A$ be a homogeneous smooth function of degree $m$. Assume that
$-\frac{d}{2}< s_1, t_1, s_2, t_2\le1+\frac{d}{2}$, Then there hold:
if\, $k\ge1$,
\begin{align}\label{2.16}
|&(A(D)\Delta_k(v\cdot\nabla f),A(D)\Delta_kf)|\nonumber\\ &
\lesssim
\alpha_k2^{-k(s_2-m)}\|v\|_{B^{\frac{d}{2}+1}}\|f\|_{\widetilde{B}^{s_1, s_2}}
\|A(D)\Delta_kf\|_2,\qquad\end{align}
and if\, $k\le0$,
\begin{align}\label{2.17}
|&(A(D)\Delta_k(v\cdot\nabla f),A(D)\Delta_kf)|\nonumber\\& \lesssim
\alpha_k2^{-k(s_1-m)}\|v\|_{B^{\frac{d}{2}+1}}\|f\|_{\widetilde{B}^{s_1, s_2}}
\|A(D)\Delta_kf\|_2,\qquad\end{align}
and if\, $k\ge1$,
\begin{align}\label{2.18}
|&(A(D)\Delta_k(v\cdot\nabla f),\Delta_kg)+(\Delta_k(v\cdot\nabla g),A(D)\Delta_kf)|
\nonumber\\&\lesssim
\alpha_k\|v\|_{B^{\frac{d}{2}+1}}(2^{-kt_2}\|g\|_{\widetilde{B}^{t_1, t_2}}\|A(D)\Delta_kf\|_2
+2^{-k(s_2-m)}\|f\|_{\widetilde{B}^{s_1, s_2}}\|\Delta_kg\|_2),\qquad
\end{align}
and if\, $k\le0$,
\begin{align}\label{2.19}
|&(A(D)\Delta_k(v\cdot\nabla f),\Delta_kg)+(\Delta_k(v\cdot\nabla g),A(D)\Delta_kf)|
\nonumber\\&\lesssim
\alpha_k\|v\|_{B^{\frac{d}{2}+1}}(2^{-kt_1}\|g\|_{\widetilde{B}^{t_1, t_2}}\|A(D)\Delta_kf\|_2
+2^{-k(s_1-m)}\|f\|_{\widetilde{B}^{s_1, s_2}}\|\Delta_kg\|_2),\qquad
\end{align}
where $\displaystyle\sum_{k\in\Z}\alpha_k\le1$.
\end{Proposition}
For the proof we refer to \cite{D2}.

\setcounter{equation}{0}
\section{Existence}

In this section, we prove the existence of the solution for the 2D viscous
shallow water equations. Without loss of generality, we assume that $\bar{h}_0=1$ and $\nu=1$. Replacing
$h$ by $h+1$ in (\ref{1.1}), we  rewrite (\ref{1.1}) as
\begin{equation}\label{3.1}
\left\{
\begin{array}{ll}
h_t+\dv u+\dv(hu)=0,\\
u_t-(\na\cdot D(u)+\na \textrm{div}\,u)+u\cdot\na u-\frac{\na h}{1+h}(D(u)+\textrm{div}\, u)+\na h=0, \\
u(0,\cdot)=u_0,  \quad   h(0,\cdot)=h_0.
\end{array}
\right.
\end{equation}

\subsection{The linearized system}

In this subsection, we consider the linearized system of (\ref{3.1}):
\begin{equation}\label{3.2}
\left\{
\begin{array}{ll}
h_t+v\cdot\na h+\dv u={\cal H},\\
u_t-(\na\cdot D(u)+\na \textrm{div}\,u)+v\cdot\na u+\na h={\cal G}, \\
u(0,\cdot)=u_0,  \quad   h(0,\cdot)=h_0.
\end{array}
\right.
\end{equation}

Let us first introduce some definitions. Set
$$e_k^r(t)\triangleq(1-e^{-cr2^{2k}t})^\f1r,\quad
\omega_k(t)=\sum_{\widetilde{k}\ge k}
2^{-(\widetilde{k}-k)}(e_{\widetilde{k}}^1(t)+e_{\widetilde{k}}^2(t)),
$$
where $c$ is a positive constant which will be determined later. We remark that
$$
\omega_k(t)\le C, \quad \textrm{for}\,\, \textrm{any}\,\, k\in \Z,
$$
which will be constantly used in the following.

\begin{Def}\label{Def3.1}
Let $s\in\R$ and $T>0$. The function space $E^{s}_T$
is  defined by
$$E^{s}_T=\{f\in {\cal Z'}((0,T)\times\R^d):\,\|f\|_{E^{s}_T}<+\infty\},$$
where
\beno
\|f\|_{E^{s}_T}\triangleq\sum_{k\in Z}2^{ks}\omega_k(T)\|\Delta_k f\|_{L^\infty_T(L^2)}.
\eeno
\end{Def}

\begin{Def}{\label{Def3.2}}
Let $s_1, s_2\in\R$ and $T>0$. The function space $\widetilde{E}^{s_1,s_2}_T$
is  defined by
$$\widetilde{E}^{s_1,s_2}_T=\{f\in {\cal Z'}((0,T)\times\R^d):\,\|f\|_{\widetilde{E}^{s_1,s_2}_T}<+\infty\},$$
where
$$\|f\|_{\widetilde{E}^{s_1,s_2}_T}\triangleq \sum_{k\le0}2^{ks_1}\omega_k(T)\|\Delta_k f\|_{L^\infty_T(L^2)}
+\sum_{k\ge 1}2^{ks_2}\omega_k(T)\|\Delta_k f\|_{L^\infty_T(L^2)}.$$
\end{Def}
\begin{rmk}
If $s_1\le s_2$, then $\widetilde{E}^{s_1,s_2}_T=E^{s_1}_T\cap E^{s_2}_T$.
Otherwise, $\widetilde{E}^{s_1,s_2}_T=E^{s_1}_T+E^{s_2}_T$.
\end{rmk}

Let $(u,h)$ be a smooth solution of (\ref{3.2}). We want to establish the following {\it a-priori} estimates for $(h,u)$:
\ben\label{3.3}
&&\|u\|_{L^1_T(B^2)}+\|u\|_{L^2_T(B^1)}+\|h\|_{\widetilde{E}^{0,1}_T}\nonumber\\
&&\quad\le C\sum_{k\in \Z}\omega_k(T)E_k(0)
+C\sum_{k\in \Z}\omega_k(T)\|\Delta_k{\cal G}(t)\|_{L^1_T(L^2)}\nonumber\\
&&\qquad+C\sum_{k\ge 1}\omega_k(T)\|\nabla \Delta_k{\cal H}(t)\|_{L^1_T(L^2)}
+C\sum_{k<1}\omega_k(T)\|\Delta_k{\cal H}(t)\|_{L^1_T(L^2)}\nonumber\\
&&\qquad+C\|u\|_{L^2_T(B^1)}\|v\|_{L^2_T(B^1)}+C\|h\|_{\widetilde{E}^{0,1}_T}\|v\|_{L^1_T(B^2)},
\een
and
\ben\label{3.4}
&&\|u\|_{L^\infty_T(B^0)}+\|h\|_{L^\infty_T(\widetilde{B}^{0,1})}+\|h\|_{L^1_T(\widetilde{B}^{2,1})}\nonumber\\
&&\quad\le E_0+C\Bigl(\|{\cal H}\|_{L^1_T(\widetilde{B}^{0,1})}
+\|{\cal G}\|_{L^1_T(B^0)}+\int_0^TV'(t)(\|u(t)\|_{B^0}+\|h(t)\|_{\widetilde{B}^{0,1}})dt\Bigr),
\een
where $V(t)=\|v(t')\|_{L^1_t(B^2)}$ and
$$
E_0=\sum_{k\in\Z}E_k(0),\quad E_{k}(t)=\bigg\{
\begin{array}{ll}E_{hk}(t)
&\quad k\ge 1,\\
E_{lk}(t)&\quad k< 1,
\end{array}\bigg.
$$
with
\beno
&& E_{hk}^2(t)=\f12\|u_k(t)\|^2_{2}+\|
\nabla h_k(t)\|^2_{2}+(u_k(t), \nabla h_k(t)), \quad \textrm{and}\\
&& E_{lk}^2(t)= \f12\|u_k(t)\|^2_{2}+\f12\|
h_k(t)\|^2_{2}+\f18(u_k(t),\nabla h_k(t)).
\eeno

Let us begin with the proof of (\ref{3.3}) and (\ref{3.4}). Set
$$
u_k=\Delta_k u,\quad h_k=\Delta_k h,\quad {\cal H}^k=\Delta_k {\cal
H},\quad {\cal G}^k=\Delta_k{\cal G}.
$$
Then we get by applying the operator $\Delta_k$ to (\ref{3.2}) that
\begin{equation}\label{3.5}
\left\{
\begin{array}{ll}
\p_th_k+\Delta_k(v\cdot\na h)+\dv u_k={\cal H}_k,\\
\p_tu_k-(\na\cdot D(u_k)+\na \textrm{div}\,u_k)+\Delta_k(v\cdot\na u)+\na h_k={\cal G}_k, \\
u_k(0,\cdot)=\Delta_ku_0,  \quad   h_k(0,\cdot)=\Delta_kh_0.
\end{array}
\right.
\end{equation}
Multiplying the second equation of (\ref{3.5}) by $u_k$, and
integrating the resulting equation over $\R^2$, we obtain
\ben\label{3.6} \frac{1}{2}\frac{d}{dt}\|u_k\|^2_{2} +\f12\|\nabla
u_k\|^2_{2}+\f32\|\dv u_k\|_2^2+(\na h_k,u_k)=({\cal G}_k, u_k)
-(\Delta_k(v\cdot\nabla u), u_k). \een

In the following, we will deal with the high frequency and the low frequency of $h$ in a different manner.

\vspace{.15cm}

{\bf High  frequencies}:\, $k\ge 1$.\vspace{.15cm}

Firstly, applying $\nabla$ to the first equation of $(\ref{3.5})$,
and  multiplying it by $\na h_k$, then integrating the resulting
equation over $\R^2$, we obtain \ben\label{3.7}
\frac{1}{2}\frac{d}{dt}\|\nabla h_k\|^2_{2}+ (\nabla\dv u_k, \nabla
h_k)=(\nabla {\cal H}_{k}, \nabla h_k)-(\nabla\Delta_k(v\cdot \nabla
h), \nabla h_k).\quad \een Secondly, applying the operator $\nabla$
to the first equation of $(\ref{3.5})$ and taking the $L^2$ product
of the resulting equation with $u_k$; then taking the $L^2$ product
of  second equation of $(\ref{3.5})$ with $\nabla h_k$, we get by
summing them up that
\ben\label{3.8} &&\frac{d}{dt}(u_k, \nabla
h_k)-\|\dv u_k\|^2_{2} -2(\na \dv u_k,\nabla h_k)+\|\nabla
h_k\|^2_{2}\nonumber\\&& \quad=(\na{\cal H}_{k}, u_k)+({\cal
G}_{k}, \nabla h_k) -(\nabla\Delta_k(v\cdot \nabla h),
u_k)-(\Delta_k(v\cdot\nabla u), \nabla h_k), \een where we used the
fact that
$$
(\na\cdot D(u_k)+\na \textrm{div}\,u_k, \na h_k)=2(\na \dv u_k,\nabla h_k).
$$
Then we get by summing up (\ref{3.6}), (\ref{3.7})$\times 2$, and (3.8) that
\ben\label{3.9}
&&\frac{d}{dt}\bigl[\f12\|u_k\|^2_{2}+\|
\nabla h_k\|^2_{2}+(u_k, \nabla h_k)\bigr]\nonumber\\
&&\qquad\qquad+\bigl[\|\nabla h_k\|^2_{2}+\f12\|\nabla u_k\|^2_{2}+\f12\|\dv u_k\|^2_2
+(\nabla h_k, u_k)\bigr]\nonumber\\
&&=\Bigl[(\nabla {\cal H}_{k}, u_k)+2(\nabla {\cal H}_{k},\na h_k)
+({\cal G}_{k}, u_k)+({\cal G}_{k}, \na h_k)\Bigr]\nonumber\\
&&\quad-(\Delta_k(v\cdot\nabla u), u_k)
-2(\nabla\Delta_k(v\cdot\nabla h), \nabla h_k)\nonumber\\
&&\quad-\Bigl[(\nabla\Delta_k(v\cdot\nabla h), u_k)
+(\Delta_k(v\cdot\nabla u), \nabla h_k)\Bigr]\nonumber\\
&&\triangleq I+II+III+IV.
\een
Note that
\beno
(u_k,\nabla h_k)\le \f13\|u_k\|^2_{2}+
\f34\|\nabla h_k\|^2_{2},
\eeno
hence, we get by the definition of $E_{hk}$ that
\ben\label{3.10}
\f16(\|u_k\|^2_{2}+\|\nabla h_k\|^2_{2})\le E_{hk}^2\le 2(\|u_k\|^2_{2}+\|\nabla h_k\|^2_{2}).
\een
Similarly, using the fact that $\f562^k\ge \f53$ and (\ref{3.10}), we have
\ben\label{3.11}
\|\nabla h_k\|^2_{2}+\f12\|\nabla u_k\|^2_{2}+\f12\|\dv u_k\|^2_{2}
+(\nabla h_k, u_k)\ge\frac{1}{8}E_{hk}^2.
\een

By summing up (\ref{3.9})-(\ref{3.11}), we obtain
\ben\label{3.12}
\frac{d}{dt}E_{hk}^2+cE_{hk}^2\le C|I+II+III+IV|.
\een

In order to obtain (\ref{3.3}), we use Lemma \ref{Lem5.1} to deal with the right hand terms of (\ref{3.12}).
Firstly, we get by using the Cauchy-Schwartz inequality and (\ref{3.10}) that
\ben\label{3.13}
|I|\le C(\|\nabla {\cal H}_{k}(t)\|_{2}
+\|{\cal G}_{k}(t)\|_{2})E_{hk}.
\een
From Lemma \ref{Lem5.1} and (\ref{3.10}), it follows that
\ben\label{3.14}
|II+III+IV|\le C(\|{\cal F}_k^1(t)\|_{2}+\|{\cal \widetilde{F}}^0_k(t)\|_{2})E_{hk}.
\een

By summing up (\ref{3.12}), and (\ref{3.13})-(\ref{3.14}), we obtain
\ben\label{3.15}
\frac{d}{dt}E_{hk}+cE_{hk}\le C\Bigl(\|\nabla {\cal H}_{k}(t)\|_{2}
+\|{\cal G}_{k}(t)\|_{2}+\|{\cal F}_k^1(t)\|_{2}+\|{\cal \widetilde{F}}^0_k(t)\|_{2}\Bigr),
\een
which implies that
\begin{align}\label{3.16}
\|E_{hk}(t)\|_{L^\infty_T}\le & E_{hk}(0)+C\Bigl(\|\nabla {\cal H}_{k}(t)\|_{L^1_T(L^2)}
+\|{\cal G}_{k}(t)\|_{L^1_T(L^2)}\nonumber\\
&+\|{\cal F}_k^1(t)\|_{L^1_T(L^2)}+\|{\cal \widetilde{F}}_k^0(t)\|_{L^1_T(L^2)}\Bigr).
\end{align}
Furthermore, by (\ref{5.1}) and (\ref{5.2}), there holds
\beno
\sum_{k\in \Z}\omega_k(T)\Big(\|{\cal F}_k^1(t)\|_{L^1_T(L^2)}+\|{\cal \widetilde{F}}_k^0(t)\|_{L^1_T(L^2)}\Big)\le C\Big(\|u\|_{L^2_T(B^1)}\|v\|_{L^2_T(B^1)}
+\|h\|_{E^{1}_T}\|v\|_{L^1_T(B^2)}\Big).
\eeno
Multiplying $\omega_k(T)$ on both sides of \eqref{3.16}, then
summing up  the resulting equation over $k\ge 1$, we obtain
\begin{align}\label{3.17}
\sum_{k\ge 1}\omega_k(T)\|E_{hk}(t)\|_{L^\infty_T}\le &\sum_{k\ge 1}\omega_k(T)E_{hk}(0)\nonumber\\
&+C\sum_{k\ge 1}\omega_k(T)\Big(\|\nabla {\cal H}_{k}(t)\|_{L^1_T(L^2)}
+\|{\cal G}_{k}(t)\|_{L^1_T(L^2)}\Big)\nonumber\\
&+C\Big(\|u\|_{L^2_T(B^1)}\|v\|_{L^2_T(B^1)}+\|h\|_{E^{1}_T}\|v\|_{L^1_T(B^2)}\Big).
\end{align}

Next, we use the decay effect of the parabolic operators to  estimate  $\|u\|_{L^2_T(B^1)\cap L^1_T(B^2)}$. It follows from (\ref{3.6}) and Lemma 5.1 that
\beno
\frac{d}{dt}\|u_k\|_2+c2^{2k}\|u_k\|_2\le C(\|\na h_k(t)\|_2+\|{\cal G}_{k}(t)\|_{L^2}
+\|{\cal \widetilde{F}}_{k}^0(t)\|_{L^2}),
\eeno
which implies that
\beno
\|u_k\|_2\le e^{-ct2^{2k}}\|u_k(0)\|_{2}+Ce^{-ct2^{2k}}\ast_t\Bigl(\|\na h_k(t)\|_2
+\|{\cal G}_{k}(t)\|_{2}+\|{\cal \widetilde{F}}_{k}^0(t)\|_{2}\Bigr),
\eeno
where the sign $\ast$ denotes the convolution of functions defined in $\R^+$, more precisely,
$$e^{-ct2^{2k}}\ast_t f\triangleq \int_0^t
e^{-c(t-\tau)2^{2k}}f(\tau)d\tau.$$
Taking the $L^r$ norm for $r=1,2$ with respect to $t$, we get by using the Young's inequality that
\beno
\|u_k\|_{L^r_T(L^2)}\le C2^{-2k/r}e^r_k(T)\Big(\|u_k(0)\|_{2}+
\|\na h_k\|_{L^1_T(L^2)}+\|{\cal G}_{k}\|_{L^1_T(L^2)}+\|{\cal \widetilde{F}}_{k}^0\|_{L^1_T(L^2)}\Big),
\eeno
which together with (\ref{5.2}) implies that
\begin{align}\label{3.18}
&\sum_{k\ge1}\Bigl(2^{2k}\|u_k\|_{L^1_T(L^2)}+2^{k}\|u_k\|_{L^2_T(L^2)}\Bigr)\le C\sum_{k\ge1}\omega_k(T)\|u_k(0)\|_2
\nonumber\\&\qquad+C\sum_{k\ge1}\omega_k(T)\Bigl(\|\na h_k\|_{L^1_T(L^2)}+\|{\cal G}_{k}\|_{L^1_T(L^2)}\Bigr)
+C\|u\|_{L^2_T(B^1)}\|v\|_{L^2_T(B^1)},
\end{align}
where we used the fact that
\beno
e_k^1(T)+e_k^2(T)\le \omega_k(T).
\eeno
On the other hand, it follows from (\ref{3.15}) that
\beno
\|E_{hk}\|_2\le e^{-ct}E_{hk}(0)+Ce^{-ct}\ast_t\Bigl(\|\na {\cal H}_k(t)\|_2
+\|{\cal G}_{k}(t)\|_{2}+\|{\cal{F}}_{k}^1(t)\|_{2}+\|{\cal \widetilde{F}}_{k}^0(t)\|_{2}\Bigr).
\eeno
Taking the $L^1$ norm with respect to $t$, we get by using the Young's inequality that
\begin{align}\label{3.19}
\|E_{hk}\|_{L_T^1}&\le C(1-e^{-cT})E_{hk}(0)+C(1-e^{-cT})\Bigl(\|\na {\cal H}_k(t)\|_{L_T^1(L^2)}
+\|{\cal G}_{k}(t)\|_{L_T^1(L^2)}\nonumber\\
&\quad+\|{\cal{F}}_{k}^1(t)\|_{L_T^1(L^2)}+\|{\cal \widetilde{F}}_{k}^0(t)\|_{L_T^1(L^2)}\Bigr).
\end{align}
Note that for $k\ge 1$
$$
1-e^{-ct}\le 1-e^{-ct2^{2k}}\le \omega_k(t),
$$
which together with (\ref{3.19}) and Lemma \ref{Lem5.1} gives
\begin{align}\label{3.20}
\sum_{k\ge 1}\|E_{hk}\|_{L_T^1}\le& C\sum_{k\ge 1}\omega_k(T)E_{hk}(0)
+C\sum_{k\ge 1}\omega_k(T)\Bigl(\|\na {\cal H}_k(t)\|_{L_T^1(L^2)}
+\|{\cal G}_{k}(t)\|_{L_T^1(L^2)}\Bigr)\nonumber\\
&+C\Big(\|u\|_{L^2_T(B^1)}\|v\|_{L^2_T(B^1)}+\|h\|_{E^1_T}\|v\|_{L^1_T(B^2)}\Big).
\end{align}
Plugging \eqref{3.20} into (\ref{3.18}), we obtain
\begin{align}\label{3.21}
&\sum_{k\ge 1}\Bigl(2^{2k}\|u_k\|_{L^1_T(L^2)}+2^{k}\|u_k\|_{L^2_T(L^2)}\Bigr)\nonumber\\&\quad\le C\sum_{k\ge1}
\omega_k(T)E_{hk}(0)
+C\sum_{k\ge 1}\omega_k(T)\Bigl(\|\na {\cal H}_k(t)\|_{L_T^1(L^2)}
+\|{\cal G}_{k}(t)\|_{L_T^1(L^2)}\Bigr)\nonumber\\
&\qquad+C\Big(\|u\|_{L^2_T(B^1)}\|v\|_{L^2_T(B^1)}+\|h\|_{E^{1}_T}\|v\|_{L^1_T(B^2)}\Big).
\end{align}

On the other hand, in order to obtain (\ref{3.4}), we use Proposition \ref{Prop2.4} to deal with the right hand terms of (\ref{3.12}).
Applying (\ref{2.16}) with $s_1=s_2=0$ to $II$, (\ref{2.16}) with $s_1=0, s_2=1$ to $III$, (\ref{2.18})
with $t_1=t_2=0, s_1=0, s_2=1$ to $IV$, we obtain
\ben\label{3.22}
|II+III+IV|\le CE_{hk}\alpha_kV'(t)(\|u\|_{B^0}+\|h\|_{\widetilde{B}^{0,1}}),
\een
with $\displaystyle\sum_{k\in \Z}\alpha_k\le 1$ and $V(t)=\|v(t')\|_{L^1_t(B^2)}$. From (\ref{3.13}) and (\ref{3.22}), it follows that
\beno
\frac{d}{dt}E_{hk}+cE_{hk}\le C\Bigl(\|\nabla {\cal H}_{k}(t)\|_{2}
+\|{\cal G}_{k}(t)\|_{2}+\alpha_kV'(t)(\|u\|_{B^0}+\|h\|_{\widetilde{B}^{0,1}})\Bigr),
\eeno
from which, a similar proof of (\ref{3.21}) ensures that
\ben\label{3.23}
&&\sum_{k\ge 1}\Bigl(\|E_{hk}\|_{L_T^1}+\|E_{hk}\|_{L_T^\infty}\Bigr)\le C\sum_{k\ge 1}E_{hk}(0)\nonumber\\
&&\qquad+C\Bigl(\|{\cal H}\|_{L^1_T(\widetilde{B}^{0,1})}
+\|{\cal G}\|_{L^1_T(B^0)}+\int_0^TV'(t)(\|u(t)\|_{B^0}+\|h(t)\|_{\widetilde{B}^{0,1}})dt\Bigr).
\een

\vspace{.15cm}

{\bf Low frequencies}:\, $k<1$.\vspace{.15cm}

Multiplying the first equation of $(\ref{3.5})$ by $h_k$, we get by integrating
the resulting equation over $\R^2$ that
\ben\label{3.24}
\frac{1}{2}\frac{d}{dt}\|h_k\|^2_{2}+
(\dv u_k, h_k)=({\cal H}_{k}, h_k)-(\Delta_k(v\cdot \nabla h),h_k).
\een
Summing up (\ref{3.6}), $(\ref{3.8})\times\frac{1}{8}$, and (\ref{3.24}), we obtain
\begin{align}\label{3.25}
&\frac{d}{dt}\Bigl[\f12\|u_k\|^2_{2}+\f12\|
h_k\|^2_{2}+\frac{1}{8}(u_k,\nabla h_k)\Bigr]
\nonumber\\
&\qquad+\Bigr[\frac{1}{8}\|\nabla h_k\|^2_{2}
+\f12\|\nabla u_k\|^2_{2}+\frac{11}{8}\|\dv u_k\|^2_{2}
-\frac{1}{4}(\na\dv\, u_k,\nabla h_k)\Bigr]\nonumber\\&=\Bigl[\f18(\nabla {\cal H}_{k}, u_k)+({\cal H}_{k},h_k)
+({\cal G}_{k}, u_k)+\f18({\cal G}_{k}, \na h_k)\Bigr]\nonumber\\
&\quad-(\Delta_k(v\cdot\nabla u), u_k)
-(\Delta_k(v\cdot\nabla h),h_k)\nonumber\\
&\quad-\f18\Bigl[(\nabla\Delta_k(v\cdot\nabla h), u_k)
+(\Delta_k(v\cdot\nabla u), \nabla h_k)\Bigr]\nonumber\\
&\triangleq I+II+III+IV.
\end{align}
Note that $2^k\le 1$, we get by the Cauchy-Schwartz inequality that
\beno
\f18(u_k,\nabla h_k)\le \f3{10}\|u_k\|_{2}\|h_k\|_{2}\le \f14\|u_k\|^2_{2}+
\f14\|h_k\|^2_{2},
\eeno
hence, we get by the definition of $E_{lk}$ that
\ben\label{3.26}
\f14(\|u_k\|^2_{2}+\|h_k\|^2_{2})\le E_{lk}^2\le 2(\|u_k\|^2_{2}+\|h_k\|^2_{2}).
\een
Similarly, we can prove
\beno
\frac{1}{4}(\na\dv\, u_k,\nabla h_k)\le \f3{5}\|\na u_k\|_{2}\|\na h_k\|_{2}\le \f{9} {10}\|\na u_k\|_2^2+\f1 {10}\|\na h_k\|_2^2,
\eeno
which together with (\ref{3.26}) implies that
\ben\label{3.27}
&&\frac{1}{8}\|\nabla h_k\|^2_{2}
+\f12\|\nabla u_k\|^2_{2}+\frac{11}{8}\|\dv u_k\|^2_{2}
-\frac{1}{4}(\na\dv\, u_k,\nabla h_k)\nonumber\\
&&\qquad\ge \frac{1}{160}2^{2k}(\|u_k\|^2_{2}+\|
h_k\|^2_{2})\ge \frac{1}{320}2^{2k}E_{lk}^2.
\een

By summing up (\ref{3.25})-(\ref{3.27}), we obtain
\ben\label{3.28}
\frac{d}{dt}E_{lk}^2+c2^{2k}E_{lk}^2\le C|I+II+III+IV|.
\een

In order to obtain (\ref{3.3}), we use Lemma 5.1 to estimate the right hand terms of (\ref{3.28}).
Using the fact that $2^k\le 1$, we get by the Cauchy-Schwartz inequality and (\ref{3.26}) that
\ben\label{3.29}
|I|\le C(\|{\cal H}_{k}(t)\|_{2}
+\|{\cal G}_{k}(t)\|_{2})E_{lk}.
\een
Using Lemma 5.1 and (\ref{3.26}), we have
\ben\label{3.30}
|II+III+IV|\le C(\|{\cal F}_k^1(t)\|_{2}+\|{\cal F}^0_k(t)\|_{2}+\|{\cal \widetilde{F}}^0_k(t)\|_{2})E_{lk}.
\een

By summing up (\ref{3.28})-(\ref{3.30}), we obtain
\beno
\frac{d}{dt}E_{lk}+c2^{2k}E_{lk}\le C\Bigl(\|{\cal H}_{k}(t)\|_{2}
+\|{\cal G}_{k}(t)\|_{2}+\|{\cal F}_k^1(t)\|_{2}+\|{\cal F}^0_k(t)\|_{2}+\|{\cal \widetilde{F}}^0_k(t)\|_{2}\Bigr),
\eeno
which implies that
\beno
E_{lk}\le e^{-c2^{2k}t}E_{lk}(0)+Ce^{-c2^{2k}t}\ast_t\Bigl(\|{\cal H}_{k}(t)\|_{2}
+\|{\cal G}_{k}(t)\|_{2}+\|{\cal F}_k^1(t)\|_{2}+\|{\cal F}^0_k(t)\|_{2}+\|{\cal \widetilde{F}}^0_k(t)\|_{2}\Bigr).
\eeno
Taking the $L^r$ norm with respect to $t$, we get by using the Young's inequality that
\begin{align}
\|E_{lk}\|_{L^r_T}&\le C2^{-2k/r}e^r_k(T)\Bigl(E_{lk}(0)+\|{\cal H}_{k}(t)\|_{L^1_T(L^2)}
+\|{\cal G}_{k}(t)\|_{L^1_T(L^2)}\nonumber\\
&\quad+\|{\cal F}_k^1(t)\|_{L^1_T(L^2)}+\|{\cal F}^0_k(t)\|_{2}+\|{\cal \widetilde{F}}^0_k(t)\|_{L^1_T(L^2)}\Bigr),\nonumber
\end{align}
from which and Lemma 5.1, it follows that
\ben\label{3.31}
&&\sum_{k<1}\omega_k(T)\|E_{lk}\|_{L^\infty_T}\le C\sum_{k<1}\omega_k(T)E_{lk}(0)+\sum_{k<1}\omega_k(T)(\|{\cal H}_{k}(t)\|_{L^1_T(L^2)}
+\|{\cal G}_{k}(t)\|_{L^1_T(L^2)})\nonumber\\
&&\qquad\quad +C\Big(\|u\|_{L^2_T(B^1)}\|v\|_{L^2_T(B^1)}+\|h\|_{\widetilde{E}^{0,1}_T}\|v\|_{L^1_T(B^2)}\Big),
\een
and
\begin{align}\label{3.32}
&\sum_{k<1}(2^{2k}\|E_{lk}\|_{L^1_T}+2^{k}\|E_{lk}\|_{L^2_T})\nonumber\\&\quad\le C\sum_{k<1}\omega_k(T)E_{lk}(0)
+\sum_{k<1}\omega_k(T)(\|{\cal H}_{k}(t)\|_{L^1_T(L^2)}
+\|{\cal G}_{k}(t)\|_{L^1_T(L^2)})\nonumber\\
&\qquad+C\Bigl(\|u\|_{L^2_T(B^1)}\|v\|_{L^2_T(B^1)}+\|h\|_{\widetilde{E}^{0,1}_T}\|v\|_{L^1_T(B^2)}\Bigr).
\end{align}

On the other hand, in order to obtain (\ref{3.4}),
we use Proposition \ref{Prop2.4} to deal with the right hand terms of (\ref{3.28}).
Applying (\ref{2.17}) with $s_1=s_2=0$ to $II$, (\ref{2.17}) with $s_1=0, s_2=1$ to $III$, (\ref{2.19})
with $t_1=t_2=0, s_1=0, s_2=1$ to $IV$, we obtain
\begin{align}\label{3.33}
|II+III+IV|\le CE_{lk}\alpha_kV'(t)(\|u\|_{B^0}+\|h\|_{\widetilde{B}^{0,1}}),
\end{align}
with $\displaystyle\sum_{k\in \Z}\alpha_k\le 1$ and $V(t)=\|v(t')\|_{L^1_t(B^2)}$. From (\ref{3.32}) and (\ref{3.36}), it follows that
\beno
\frac{d}{dt}E_{lk}+c2^{2k}E_{lk}\le C\Bigl(\|{\cal H}_{k}(t)\|_{2}
+\|{\cal G}_{k}(t)\|_{2}+\alpha_kV'(t)(\|u\|_{B^0}+\|h\|_{\widetilde{B}^{0,1}})\Bigr),
\eeno
from which and a similar proof of (\ref{3.21}) ensure that
\ben\label{3.34}
&&\sum_{k<1}\Bigl(2^{2k}\|E_{lk}\|_{L_T^1}+\|E_{lk}\|_{L_T^\infty}\Bigr)\le \sum_{k<1}E_{hk}(0)\nonumber\\
&&\qquad+C\Bigl(\|{\cal H}\|_{L^1_T(\widetilde{B}^{0,1})}
+\|{\cal G}\|_{L^1_T(B^0)}+\int_0^TV'(t)(\|u(t)\|_{B^0}+\|h(t)\|_{\widetilde{B}^{0,1}})dt\Bigr).
\een

\vspace{.2cm}
\noindent{\bf The completion of the {\it a-priori} estimates}\vspace{.2cm}

\no Firstly, adding up \eqref{3.17}, \eqref{3.21}, \eqref{3.31}, and \eqref{3.32} yields that
\ben\label{3.35}
&&\|u\|_{L^1_T(B^2)}+\|u\|_{L^2_T(B^1)}+\|h\|_{\widetilde{E}^{0,1}_T}\nonumber\\
&&\quad\le C\sum_{k\in \Z}\omega_k(T)E_k(0)
+C\sum_{k\in \Z}\omega_k(T)\|{\cal G}_k(t)\|_{L^1_T(L^2)}\nonumber\\
&&\qquad+C\sum_{k\ge 1}\omega_k(T)\|\nabla{\cal H}_k(t)\|_{L^1_T(L^2)}
+C\sum_{k<1}\omega_k(T)\|{\cal H}_k(t)\|_{L^1_T(L^2)}\nonumber\\
&&\qquad+C\|u\|_{L^2_T(B^1)}\|v\|_{L^2_T(B^1)}+C\|h\|_{\widetilde{E}^{0,1}_T}\|v\|_{L^1_T(B^2)},
\een
where we used the fact that
$$
\|h\|_{E^{1}_T}\le C\|h\|_{\widetilde{E}^{0,1}_T}.
$$
On the other hand, adding up \eqref{3.23} and \eqref{3.34} gives rise to
\ben\label{3.36}
&&\|u\|_{L^\infty_T(B^0)}+\|h\|_{L^\infty_T(\widetilde{B}^{0,1})}+\|h\|_{L^1_T(\widetilde{B}^{2,1})}
\nonumber\\&&\qquad \le  E_0+C\Big(\|{\cal H}\|_{L^1_T(\widetilde{B}^{0, 1})}+\|{\cal G}\|_{L^1_T(B^0)}+
\int_0^TV'(t)(\|u\|_{B^{0}}+\|h\|_{\widetilde{B}^{0, 1}})dt\Big),
\een
which together with the Gronwall inequality implies that
\ben\label{3.37}
&&\|u\|_{L^\infty_T(B^0)}+\|h\|_{L^\infty_T(\widetilde{B}^{0,1})}+\|h\|_{L^1_T(\widetilde{B}^{2,1})}\nonumber\\
&&\qquad\qquad\le Ce^{C\|v\|_{L^1_T(B^2)}}\Big(E_0+\|{\cal H}\|_{L^1_T(\widetilde{B}^{0, 1})}+\|{\cal G}\|_{L^1_T(B^0)}\Big).
\een
Finally, let us remark that
\beno
E_0\thickapprox(\|h_0\|_{\widetilde{B}^{0,1}}+\|u_0\|_{B^0}).
\eeno

\subsection{The uniform estimate of the approximate sequence of solutions}

In this subsection, we will construct the approximate solutions of (\ref{3.1})
and present the uniform estimate of the approximate solutions.
Let us first define the approximate sequence $(h^n, u^n)_{n\in\N}$ of (\ref{3.1}) by the following system:
\begin{align}\label{3.38} \left\{
\begin{aligned}
&\partial_t h^{n+1}+u^n\cdot\nabla h^{n+1}+\dv u^{n+1}=
{\cal H}^n,\\
&\partial_tu^{n+1}-(\na\cdot D(u^{n+1})+\na \textrm{div}\,u^{n+1})+u^{n}\cdot\nabla u^{n+1}+\nabla h^{n+1}
={\cal G}^{n},\qquad\quad\\
&(h^{n+1}, u^{n+1})|_{t=0}=\sum_{|k|\le n+N}\Delta_k(h_0,u_0),
\end{aligned}
\right.
\end{align}
where
$${\cal H}^n\triangleq -h^n\dv u^n,\qquad {\cal G}^{n}\triangleq \frac{\nabla h^n}{1+h^n}\widetilde{\nabla }u^n,
\quad \textrm{with}\quad \widetilde{\na} u^n=D(u^n)+\dv u^n, $$
and $N$ is a fixed large integer such that
$$
1+h^n(0)\ge \f34,\qquad \textrm{for} \quad n\ge 1.
$$
Set $(h^0, u^0)=(0,0)$ and solve the linear system, we can define  $(h^n, u^n)_{n\in\N_0}$ by the induction.
Next, we are going to prove by the induction that there exist positive constants $\eta$, $K$, and  $T$ such that
the following bounds hold for all $n\in\N_0$:
\begin{align}
&1+h^n\ge\f12,\label{3.39}\\
&\|u^{n}\|_{L^1_T(B^2)\cap L^2_T(B^1)}+\|h^n\|_{\widetilde{E}^{0,1}_T}
\le \eta,\label{3.40}\\
&\|u^n\|_{L^\infty_T(B^0)}+\|h^n\|_{L^\infty_T(\widetilde{B}^{0,1})\cap L^1_T(\widetilde{B}^{2,1})}
\le KE_0.\label{3.41}
\end{align}
Assume that (\ref{3.39})-(\ref{3.41}) hold for $(h^n, u^n)$,
we need to prove that (\ref{3.39})-(\ref{3.41}) also hold for $(h^{n+1}, u^{n+1})$.
Applying the {\it a-priori} estimates (\ref{3.35}) and (\ref{3.37}) to $(h^{n+1}, u^{n+1})$, we obtain
\ben\label{3.42}
&&\|u^{n+1}\|_{L^1_T(B^2)}+\|u^{n+1}\|_{L^2_T(B^1)}+\|h^{n+1}\|_{\widetilde{E}^{0,1}_T}\nonumber\\
&&\quad\le C{\cal Q}_0(T)+C\sum_{k\in \Z}\omega_k(T)\|{\cal G}_k^{n}(t)\|_{L^1_T(L^2)}
+C\sum_{k\ge 1}\omega_k(T)\|\nabla {\cal H}_k^{n}(t)\|_{L^1_T(L^2)}\nonumber\\
&&\qquad+C\sum_{k<1}\omega_k(T)\|{\cal H}_k^{n}(t)\|_{L^1_T(L^2)}+C\|u^{n+1}\|_{L^2_T(B^1)}\|u^n\|_{L^2_T(B^1)}\nonumber\\
&&\qquad+C\|h^{n+1}\|_{\widetilde{E}^{0,1}_T}\|u^n\|_{L^1_T(B^2)},
\een
and
\ben\label{3.43}
&&\|u^{n+1}\|_{L^\infty_T(B^0)}+\|h^{n+1}\|_{L^\infty_T(\widetilde{B}^{0,1})}+\|h^{n+1}\|_{L^1_T(\widetilde{B}^{2,1})}\nonumber\\
&&\qquad\qquad\le Ce^{C\|u^{n}\|_{L^1_T(B^2)}}\Big(E_0+\|{\cal H}^{n}\|_{L^1_T(\widetilde{B}^{0,1})}+\|{\cal G}^{n}\|_{L^1_T(B^0)}\Big),
\een
with
$${\cal Q}_0(T)\triangleq\sum_{k\in\Z}\omega_k(T)E_{k}(0).$$
Thanks to (\ref{2.4}), we have
\beno
\|{\cal H}^n\|_{{B}^{0}}\le C\|{h}^n\|_{{B}^{0}}
\|u^n\|_{B^2},\quad \hbox{and}\quad
\|{\cal H}^n\|_{{B}^{1}}\le C\|{h}^n\|_{{B}^{1}}
\|u^n\|_{B^2},
\eeno which together with the fact that $\widetilde{B}^{0,1}={B}^{0}\cap{B}^{1}$ yields
\ben\label{3.44}
\|{\cal H}^n\|_{L^1_T(\widetilde{B}^{0,1})}
\le C\|{h}^n\|_{L^\infty_T(\widetilde{B}^{0,1})}\|u^n\|_{L^1_T(B^2)}\le CKE_0\eta.
\een
We rewrite ${\cal G}^n$ as
\beno\label{3.45}
\frac{\nabla h^n}{1+h^n}\widetilde{\na} u^n=
(1+h^n)\nabla\bigg(\frac{h^n}{1+h^n}\bigg)\widetilde{\na} u^n.
\eeno
Using (\ref{2.4}) and  (\ref{2.12}), we get
\begin{align}\label{3.45}
\|{\cal G}^{n}\|_{L^1_T(B^0)}
&\le C\Bigl\|\nabla\bigg(\frac{h^n}{1+h^n}\bigg)\Bigr\|_{L^\infty_T(B^0)}\|(1+h^n)\widetilde{\na} u^n\|_{L^1_T(B^1)}\nonumber\\
&\le C(1+\|h^n\|_{L^\infty_T(L^\infty)})^2\|h^n\|_{L^\infty_T(B^1)}(1+\|h^n\|_{L^\infty_T(B^1)})\|u^n\|_{L^1_T(B^2)}\nonumber\\
&\le C(1+\|h^n\|_{L^\infty_T(\widetilde{B}^{0,1})})^3
\|h^{n}\|_{L^\infty_T(B^1)}\|u^n\|_{L^1_T(B^2)}\nonumber\\
&\le CKE_0(1+KE_0)^3\eta.
\end{align}
Plugging  (\ref{3.44}) and (\ref{3.45}) into (\ref{3.43}) yields that
\ben\label{3.46}
\|u^{n+1}\|_{L^\infty_T(B^0)}+\|h^{n+1}\|_{L^\infty_T(\widetilde{B}^{0,1})}+\|h^{n+1}\|_{L^1_T(\widetilde{B}^{2,1})}
\le Ce^{C\eta}\Big(E_0+KE_0(1+KE_0)^3\eta\Big).
\een
We take $T, \eta>0$ small enough and $K=4C$ such that
\begin{align}
e^{C\eta}\le 2,\quad  K(1+KE_0)^3\eta\le 1,\tag{$\Re_1$}
\end{align}
from which and (\ref{3.46}), it follows that
\beno
\|u^{n+1}\|_{L^\infty_T(B^0)}+\|h^{n+1}\|_{L^\infty_T(\widetilde{B}^{0,1})}+\|h^{n+1}\|_{L^1_T(\widetilde{B}^{2,1})}\le KE_0.
\eeno
This proves (\ref{3.41}) for $(u^{n+1},h^{n+1})$.

Next, we prove (\ref{3.40}) for $(u^{n+1},h^{n+1})$.
Applying Lemma \ref{Lem5.2} with $s_1=0$ and $s_2=1$, (\ref{2.4}) with $s_1=s_2=1$,
and Lemma \ref{Lem5.4} with $s=1$,  we obtain
\begin{align}\label{3.47}
\sum_{k\in \Z}\omega_k(T)\|{\cal G}_k^{n}(t)\|_{L^1_T(L^2)}
&\le C\Bigl\|\nabla\bigg(\frac{h^n}{1+h^n}\bigg)\Bigr\|_{E^0_T}\|(1+h^n)\widetilde{\na} u^n\|_{L^1_T(B^1)}\nonumber\\
&\le C(1+\|h^n\|_{L^\infty_T(L^\infty)})^3\|h^n\|_{E^1_T}(1+\|h^n\|_{L^\infty_T(B^1)})\|u^n\|_{L^1_T(B^2)}\nonumber\\
&\le C(1+\|h^n\|_{L^\infty_T(\widetilde{B}^{0,1})})^4
\|h^{n}\|_{\widetilde{E}^{0,1}_T}\|u^n\|_{L^1_T(B^2)}\nonumber\\
&\le C(1+KE_0)^4\eta^2.
\end{align}
On the other hand, we apply Lemma 5.2 with $s_1=0, s_2=1$ to get
\beno
&&\sum_{k\ge 1}\omega_k(T)\|\nabla {\cal H}_k^{n}(t)\|_{L^1_T(L^2)}+\sum_{k<1}\omega_k(T)\|{\cal H}_k^{n}(t)\|_{L^1_T(L^2)}\nonumber\\
&&\quad\le C\sum_{k\in \Z}\omega_k(T)(\|\na h^n_k\|_{L^\infty_T(L^2)}+
\|h^n_k\|_{L^\infty_T(L^2)})\|\dv u^n\|_{L^1_T(B^1)}\nonumber\\
&&\qquad+C\sum_{k\in\Z}\omega_k(T)2^{2k}\|u^n_k\|_{L^1_T(L^2)}\|h^n\|_{L^\infty_T(\widetilde{B}^{0,1})}\nonumber\\
&&\quad\triangleq I+II.
\eeno
Obviously, we have
\ben\label{3.48}
I\le C\|h^n\|_{\widetilde{E}^{0,1}_T}\|u^n\|_{L^1_T(B^2)}\le C\eta^2.
\een
In order to estimate $II$, we first fix $k_0\ge 1$ such that
\ben\label{3.49}
\sum_{k\ge k_0}\|u_k(0)\|_{2}\le \f {\eta} {16CKE_0}.
\een
Then we write
\begin{align}
II&=\sum_{k\ge k_0}\omega_k(T)2^{2k}\|u^n_k\|_{L^1_T(L^2)}\|h^n\|_{L^\infty_T(\widetilde{B}^{0,1})}
+\sum_{k\le k_0}\omega_k(T)2^{2k}\|u^n_k\|_{L^1_T(L^2)}\|h^n\|_{L^\infty_T(\widetilde{B}^{0,1})}\nonumber\\
&\triangleq II_1+II_2.\nonumber
\end{align}
Using (\ref{3.18}), (\ref{3.47}), and (\ref{3.49}), we obtain
\begin{align}\label{3.50}
II_1&\le CKE_0\Bigl[\sum_{k\ge k_0}\omega_k(T)\|u_k(0)\|_2+
\sum_{k\ge k_0}\omega_k(T)\Bigl(\|\na h_k^n\|_{L^1_T(L^2)}
+\|{\cal G}_{k}^{n-1}\|_{L^1_T(L^2)}\Bigr)\nonumber\\
&\quad+\|u^n\|_{L^2_T(B^1)}\|u^{n-1}\|_{L^2_T(B^1)}\Bigr]\nonumber\\
&\le CKE_0\Bigl[\f {\eta} {16CKE_0}+
\sum_{k\ge k_0}\omega_k(T)\|\na h_k^n\|_{L^1_T(L^2)}+(1+KE_0)^4\eta^2\Bigr].
\end{align}
On the other hand, thanks to (\ref{3.19}) and Lemma \ref{Lem5.1}, we have
\begin{align}
\sum_{k\ge k_0}\omega_k(T)\|\na h_k^n\|_{L_T^1(L^2)}&\le C(1-e^{-cT})\sum_{k\ge k_0}\omega_k(T)E_{hk}(0)\nonumber\\
&\quad+C(1-e^{-cT})\sum_{k\ge k_0}\omega_k(T)\Bigl(\|\na {\cal H}_k^{n-1}(t)\|_{L_T^1(L^2)}
+\|{\cal G}_{k}^{n-1}(t)\|_{L_T^1(L^2)}\Bigr)\nonumber\\
&\quad+C\|u^n\|_{L^2_T(B^1)}\|u^{n-1}\|_{L^2_T(B^1)}+C\|h^n\|_{E^1_T}\|u^{n-1}\|_{L^2_T(B^1)}\nonumber\\
&\le C(1-e^{-cT})E_0+C(1-e^{-cT})KE_0\eta+C(1+KE_0)^4\eta^2,\nonumber
\end{align}
where we used (\ref{3.44}) and (\ref{3.47}) in the second inequality.
Plugging the above inequality into (\ref{3.50}) yields that
\ben\label{3.51}
II_1\le CKE_0\Bigl[\f {\eta} {16CKE_0}+
(1-e^{-cT})(E_0+KE_0\eta)+(1+KE_0)^4\eta^2\Bigr].
\een
Note for $k\le k_0$, we can choose $T>0$ small enough so that
\begin{align}
\omega_k(T)\le \f 1{16CKE_0\eta},\tag{$\Re_2$}
\end{align}
so we get
\ben\label{3.52}
|II_2|\le \f \eta {16}.
\een

Plugging (\ref{3.47}), (\ref{3.48}), (\ref{3.51}), (\ref{3.52}) into (\ref{3.42}), we get
\ben\label{3.53}
&&\|u^{n+1}\|_{L^1_T(B^2)}+\|u^{n+1}\|_{L^2_T(B^1)}+\|h^{n+1}\|_{\widetilde{E}^{0,1}_T}\nonumber\\
&&\quad\le C{\cal Q}_0(T)+\f \eta {8}+ C(1+KE_0)^5\eta^2+
CKE_0(1-e^{-cT})(E_0+KE_0\eta)\nonumber\\
&&\qquad+C\eta(\|u^{n+1}\|_{L^2_T(B^1)}+\|h^{n+1}\|_{\widetilde{E}^{0,1}_T}).
\een
Note that ${\cal Q}_0(0)=0$, we can take $T, \eta$ small enough such that
\begin{align}
&C\eta\le \f12,\quad C{\cal Q}_0(T)\le \f \eta 8,\quad C(1+KE_0)^5\eta<\f18,
\quad \textrm{and}\nonumber\\
&CKE_0(1-e^{-cT})(E_0+KE_0\eta)\le \f \eta 8, \tag{$\Re_3$}
\end{align}
which together with (\ref{3.53}) gives
\beno
\|u^{n+1}\|_{L^1_T(B^2)}+\|u^{n+1}\|_{L^2_T(B^1)}+\|h^{n+1}\|_{\widetilde{E}^{0,1}_T}\le\eta.
\eeno

Finally, let us prove (\ref{3.39}) for $h^{n+1}$.
We  rewrite the first equation of \eqref{3.38} as
\beno
\partial_t(1+h^{n+1})+u^{n}\cdot\na(1+h^{n+1})+\dv u^{n+1}-{\cal H}^{n}=0.
\eeno
Then $1+h^{n+1}$ can be represented as
\begin{align}\label{3.54}
(1+h^{n+1})(t,x)=&(1+h^{n+1}_0)((\psi^n)^{-1}_t(x))+\int_0^t\dv\,u^{n+1}(\tau,\psi^n_\tau((\psi^n)^{-1}_t(x)))d\tau\nonumber\\
&+\int_0^t{\cal H}^{n}(\tau,\psi^n_\tau((\psi^n)^{-1}_t(x)))d\tau,
\end{align}
where the flow map $\psi^{n}_t$ is defined by
\beno\left\{
\begin{aligned}
&\partial_t\psi^{n}_t(x)=u^{n}(t,\psi^{n}_t(x))\\
&\psi^{n}_t|_{t=0}=x.
\end{aligned}\right.
\eeno
Thanks to the inclusion map $B^1\hookrightarrow L^\infty$ and (\ref{2.4}), we get
\beno
&&\int_0^t\|\dv\,u^{n+1}(\tau,\psi^n_\tau((\psi^n)^{-1}_t(x)))\|_{\infty} d\tau\le \|u^{n+1}\|_{L^1_t(B^2)}\le \eta,\\
&&\int_0^t\|{\cal H}^{n}(\tau,\psi_\tau(\psi^{-1}_t(x)))\|_\infty d\tau \le\|h^{n}\dv u^{n}\|_{L^1_t(B^1)}\nonumber\\
&&\qquad\qquad\qquad\qquad\le C\|h^{n}\|_{L^\infty_t(\widetilde{B}^{0,1})}\|u^{n}\|_{L^1_t(B^2)}\le CKE_0\eta,
\eeno
from which and (\ref{3.54}), it follows that
\ben\label{3.55}
1+h^{n+1}\ge \f34-(1+CKE_0)\eta.
\een
We take $\eta$ small enough such that
\begin{align}
(1+CKE_0)\eta\le \f14, \tag{$\Re_4$}
\end{align}
which together with (\ref{3.55}) ensures that
\beno
1+h^{n+1}\ge \f12.
\eeno
So far, we have show that $T$, $\eta$ can be chosen small enough such that
the assumption $(\Re_1)-(\Re_4)$ hold under which
the approximate solutions $(u^n, h^n)_{n\in\N_0}$ is uniformly bounded in
$${\cal E}_T\triangleq \Bigl(L^\infty_T(B^0)\cap L^1_T(B^2)\Bigr)\times
\Bigl(L^\infty_T(\widetilde{B}^{0,1})\cap L^1_T(\widetilde{B}^{2,1})\Bigr).$$
It should be pointed out that if $\|u_0\|_{B^0}+
\|h_0\|_{\widetilde{B}^{0,1}}$ is small enough, we can take $T=+\infty$ such that
the assumption $(\Re_1)-(\Re_4)$ hold.

\subsection {The existence of the solution}
Now let us  turn to prove the existence of the solution, and
the standard compact arguments will be used.
In the section 3.2, we have showed that the approximate solutions $(h^n,u^n)_{n\in\N}$
satisfy \eqref{3.39}-\eqref{3.41}, and without loss of generality, we can assume the following:
\begin{align}
1&+h^n\ge\f12,\label{3.56}\\
\|u^n\|_{L^\infty_T(B^0)\cap L^1_T(B^2)}&+\|h^n\|_{L^\infty_T(\widetilde{B}^{0,1})\cap L^1_T(\widetilde{B}^{2,1})}
\le KE_0.\label{3.57}
\end{align}
Using the interpolation and the fact that $B^0\cap B^1=\widetilde{B}^{0,1}$, we have
\beno
&&\|h^n\|_{L^2_T(B^1)}\lesssim\|h^n\|^{\frac{1}{2}}_{L^\infty_T(\widetilde{B}^{0,1})}
\|h^n\|^{\frac{1}{2}}_{L^1_T(\widetilde{B}^{2,1})},\quad
\|u^n\|_{L^2_T(B^1)}\lesssim
\|u^n\|^{\frac{1}{2}}_{L^\infty_T(B^0)}\|u^n\|^{\frac{1}{2}}_{L^1_T(B^2)},\\
&&\|h^n\|_{L^4_T(B^{\frac{1}{2}})}\lesssim\|h^n\|^{\frac{1}{2}}_{L^\infty_T(\widetilde{B}^{0,1})}
\|h^n\|^{\frac{1}{2}}_{L^2_T(B^{1})},\quad\|u^n\|_{L^\frac{4}{3}_T(B^\frac{3}{2})}\lesssim
\|u^n\|^{\frac{1}{4}}_{L^\infty_T(B^0)}\|u^n\|^{\frac{3}{4}}_{L^1_T(B^2)},
\eeno from which and (\ref{3.56}), it follows that
\ben\label{3.58}
\|h^n\|_{L^2_T(B^1)}+\|u^n\|_{L^2_T(B^1)}+\|h^n\|_{L^4_T(B^{\frac{1}{2}})}+\|u^n\|_{L^\frac{4}{3}_T(B^\frac{3}{2})}
\lesssim KE_0.
\een

Now, we  show that $(h^n, u^n)$ is uniformly bounded in $C^{\frac{1}{2}}_{loc}(B^0)\times
C^{\frac{1}{4}}_{loc}(B^{-\frac{1}{2}})$. Using (\ref{2.4}), (\ref{3.57})
and (\ref{3.58}), it is easy to verify that
\beno
&&\|u^{n}\cdot\nabla h^{n+1}\|_{L^2_T(B^{0})}\lesssim
\|u^{n}\|_{L^2_T(B^{1})}\|h^{n+1}\|_{L^\infty_T(\widetilde{B}^{0,1})}\lesssim (KE_0)^2,\\
&&\|h^{n}\dv u^{n}\|_{L^2_T(B^{0})}\lesssim
\|u^{n}\|_{L^2_T(B^{1})}\|h^{n}\|_{L^\infty_T(\widetilde{B}^{0,1})}\lesssim
(KE_0)^2, \eeno
from which and the first equation of (\ref{3.38}), it follows that $\partial_t
h^{n}$ is uniformly bounded in $L^2_T(B^0)$ which implies $h^{n}$ is
uniformly bounded in $C^{\frac{1}{2}}_{loc}(B^0)$.
On the other hand, thanks to (\ref{2.4}), (\ref{3.56}) and (\ref{2.12}), we have
\beno
&&\|u^{n}\cdot\nabla u^{n+1}\|_{L^\frac{4}{3}_T(B^{-\frac{1}{2}})}\lesssim
\|u^{n}\|_{L^\infty_T(B^0)}\|u^{n+1}\|_{L^\frac{4}{3}_T(B^\frac{3}{2})}\lesssim (KE_0)^2,\\
&&\bigg\|\frac{\nabla h^{n}}{1+h^{n}}\widetilde{\nabla} u^{n}\bigg\|_{L^\frac{4}{3}_T(B^{-\frac{1}{2}})}
\lesssim C(1+\|h^{n}\|_{L^\infty_T(B^1)})^3\|u^{n}\|_{L^\frac{4}{3}_T(B^\frac{3}{2})}\lesssim C(1+KE_0)^3KE_0,
\eeno
from which and the second equation of (\ref{3.38}), it follows that
$\partial_t u^{n}$ is uniformly bounded in $L^\frac{4}{3}_T(B^{-\frac{1}{2}})$
which implies $u^{n}$ is uniformly bounded
in $C^{\frac{1}{4}}_{loc}(B^{-\frac{1}{2}})$.

Next, we claim that the inclusions  $B^0\cap B^1\hookrightarrow L^2$ and
$B^{-\frac{1}{2}}\cap B^0\hookrightarrow \dot{H}^{-\frac{1}{2}}$ are locally
compact. Indeed, these can be proved by noting that for $s'<s$,
$\dot{H}^{s'}\cap \dot{H}^{s}\hookrightarrow \dot{H}^{s'}$ is locally compact and
for $s\in\R$,
$B^s\hookrightarrow \dot{H}^s$.
Then, by the Arzela-Ascoli theorem and Cantor's diagonal process, there
exist a subsequence $(u^{n_k}, h^{n_k})$ and a function $(u, h)$ such that
\ben\label{3.59}
(u^{n_k}, h^{n_k})\rightarrow (u, h)\qquad
\mbox{in}\qquad C_{loc}(\dot{H}^{-\frac{1}{2}}_{loc})\times C_{loc}(L^2_{loc}),
\een
as $n_k\rightarrow\infty$.
On the other hand, $(u^{n_k}, h^{n_k})$ is uniformly bounded
in ${\cal E}_T$, then there exists a
subsequence (which  still denoted by $(u^{n_k},h^{n_k})$)
such that
$$(u^{n_k}, h^{n_k})\rightharpoonup (u, h)\quad
\mbox{in}\quad {\cal E}_{T},$$
where $``\rightharpoonup"$ denotes weak* convergence.

Finally, let us  prove that $(u, h)$ solves (\ref{1.1})
in the sense of distribution. We only need to
prove the nonlinear terms  such as  $u^n\cdot\nabla h^n$,
$\frac{\nabla h^n}{1+h^n}\widetilde{\nabla} u^n$, etc tend
to the corresponding nonlinear terms  in the sense of distribution.
This can be done by using the uniform
estimates of $(u^n, h^n)$, $(u, h)$ in ${\cal E}_T$ and the convergence result (\ref{3.59}).
Here, we only show the case of the term
$Y(h^n)\widetilde{\nabla} u^n$
(where $Y(z)\triangleq \nabla z/(1+z) $), the other terms can be treated in the same way.
For any test function $\theta\in C_0^\infty([0, T^*)\times\R^2)$, we write
\begin{align}
\big<&Y(h^n)\widetilde{\nabla} u^n
-Y(h)\widetilde{\nabla} u,\,\theta\big>\nonumber\\
&=
\Big<(1+h^n)\nabla\Big(\frac{h^n}{1+h^n}-\frac{h}{1+h}\Big)\widetilde{\nabla} u^n,\,\theta\Big>\nonumber\\&\quad+
\Big<(h^n-h)\nabla\Big(\frac{h}{1+h}\Big)\widetilde{\nabla} u^n,\,\theta\Big>+
\Big<(1+h)\nabla\Big(\frac{h}{1+h}\Big)\widetilde{\nabla} (u^n-u),\,\theta\Big>\nonumber\\&
\triangleq{I}_1+{I}_2+{I}_3.\nonumber
\end{align}
Thanks to (\ref{2.4}) and (\ref{3.56}), we have
\begin{align}
{I}_1&\le\bigg\|\frac{\psi(h^n-h)}{(1+h^n)(1+h)}\bigg\|_2
\|\nabla((1+h^n)\widetilde{\nabla} u^n\theta)\|_2\lesssim\|\theta(h^n-h)\|_2\|(1+h^n)\widetilde{\nabla} u^n\|_{B^1}\nonumber\\
&\lesssim \|\theta(h^n-h)\|_2(1+\|h^n\|_{\widetilde{B}^{0,1}})\|u^n\|_{B^2},\nonumber
\end{align}
where $\psi \in C_0^\infty([0, T^*)\times\R^2)$,\,and $\psi=1$ on \,$\textrm{supp}\,\theta$.
For ${I}_2$, we have
\begin{align}
{I}_2&\le\|\theta(h^n-h)\|_2
\bigg\|\nabla\bigg(\frac{h}{1+h}\bigg)\widetilde{\nabla} u^n\bigg\|_2
\lesssim\|\theta(h^n-h)\|_2\|\nabla h\|_2\|\nabla u^n\|_{L^\infty}\nonumber\\
&\lesssim \|\theta(h^n-h)\|_2\|h\|_{\widetilde{B}^{0,1}}\|u^n\|_{B^2}.\nonumber
\end{align}
Using (\ref{3.56}) and the interpolation, we get
\begin{align}
{I}_3&\le\bigg\|(1+h)\nabla\bigg(\frac{h}{1+h}\bigg)\bigg\|_2
\|\widetilde{\nabla} (u^n-u)\theta\|_2\lesssim (1+\|h\|_\infty)\|\nabla h\|_2\|(u^n-u)\psi\|_{\dot{H}^1}\nonumber\\
&\lesssim (1+\|h\|_{\widetilde{B}^{0,1}})\|h\|_{B^1}\|u^n-u\|^{\frac{3}{5}}_{\dot{H}^2}
\|(u^n-u)\theta\|^{\frac{2}{5}}_{\dot{H}^{-\frac{1}{2}}}.\nonumber
\end{align}
Thus, by (\ref{3.59}), we get as $n\rightarrow 0$
$$
\big<Y(h^n)\widetilde{\nabla} u^n
-Y(h)\widetilde{\nabla} u,\,\theta\big>\longrightarrow 0.
$$
Following the argument in \cite{D1}, we can also prove that $(u,h)$ is continuous in
time with values in $B^{0}\times \widetilde{B}^{0,1}$.

\setcounter{equation}{0}

\section{Uniqueness}

In this section, we will prove the uniqueness of the solution. Firstly, let us recall some known results.
\begin{Lemma}\label{Lemma4.1}(Osgood's lemma) Let $\rho$ be a measurable positive
function and $\gamma$ a positive locally integrable function, each defined on
the domain $[t_0, t_1]$. Let $\mu: [0, \infty)\rightarrow [0,\infty)$
be a continuous nondecreasing function, with $\mu(0)=0$. Let $a\ge0$,
and assume that for all $t$ in $[t_0, t_1]$,
$$\rho(t)\le a+\int_{t_0}^{t}\gamma(\tau)\mu(\rho(\tau))d\tau.$$
If $a>0$, then
$$-{\cal M}(\rho(t))+{\cal M}(a)\le \int_{t_0}^{t}\gamma(\tau)d\tau,
\quad\mbox{where}\quad{\cal M}(x)=\int_x^1\frac{d\tau}{\mu(\tau)}.$$
If $a=0$ and ${\cal M}=\infty$, then $\rho\equiv0.$
\end{Lemma}
This Lemma can be understood as a generalization of classical Gronwall Lemma and
can be found in \cite{Ch1}.
\begin{Proposition}\label{Prop4.1}
Let $s\in (-\frac{d}{p}, 1+\frac{d}{p})$, and $1\le p,r\le+\infty$.
Let $v$ be a vector field such that $\nabla v\in L^1_T(\dot{B}^{\frac{d}{p}}_{p,r}\cap L^\infty)$.
Assume that $f_0\in \dot{B}^{s}_{p,r},$ $g\in L^1_T(\dot{B}^{s}_{p,r})$ and
$f\in L^\infty_T(\dot{B}^{s}_{p,r})\cap C([0,T]; {\cal S}')$ is the solution of
\begin{align}
\bigg\{\begin{aligned}
&\partial_t f+v\cdot \nabla f =g,\\
&f(0,x)=f_0.
\end{aligned}
\bigg.\nonumber\end{align}
Then there exists a constant $C(s,p,d)$ such that
for $t\in[0,T]$
\begin{align}\label{4.1}
\|f\|_{\widetilde{L}^\infty_t(\dot{B}^{s}_{p,r})}\le Ce^{CV(t)}\bigg(
\|f_0\|_{\dot{B}^{s}_{p,r}}+\int_0^t e^{-CV(\tau)}\|g(\tau)\|_{\dot{B}^{s}_{p,r}}d\tau\bigg),
\end{align}
where $V(t)\triangleq\int_0^t\|\nabla v(\tau)\|_{\dot{B}^{\frac{d}{p}}_{p,r}\cap L^\infty}d\tau.$
If $r<+\infty$, then $f$ belongs to $C([0,T]; \dot{B}^s_{p,r})$.
\end{Proposition}
The proof can be found in \cite{D4}.
\begin{Proposition}\label{Prop4.3}
Let $T>0$, $s\in \R$, and $1\le q,r\le+\infty$. Assume that
$u_0\in \dot{B}^{s}_{2,q},$ $g\in \widetilde{L}^1_T(\dot{B}^{s}_{2,q})$ and
$u$ is the solution of
\begin{align}
\,\bigg\{\begin{aligned}
&\partial_t u-\nu \widetilde{\Delta} u=g,\nonumber\\
&u(0,x)=u_0,
\end{aligned}
\bigg.\end{align}
where $\widetilde{\Delta} u=\na\cdot D(u)+\na \dv \,u$.
Then there exists a constant $C(s,d,\nu)$ such that
\begin{align}\label{4.2}
(r\nu)^{\frac1r}\|u\|_{\widetilde{L}^r_T(\dot{B}^{s+\frac2r}_{2,q})}\le
\Big(&\sum_{k\in\Z}\big({1-e^{-r\nu 2^{2k}T}}\big)^{\frac q r}2^{qks}\|\Delta_ku_0\|_2^q\Big)^{\frac1q}
\nonumber\\&+C\Big(\sum_{k\in\Z}\big({1-e^{-r\nu 2^{2k}T}}\big)^{\frac qr}2^{qks}\|\Delta_kg\|_{L^1_T(L^2)}^q\Big)^{\frac1q}.
\end{align}
If $q<+\infty$, then $u$ belongs to $C([0,T]; \dot{B}^{s}_{2,q})$.
\end{Proposition}
The proof is similar to the case when the diffusion term $\widetilde{\Delta} u$ is replaced by $\Delta u$.
We can refer to \cite{Ch2} see the details.

Now we introduce the logarithmic interpolation inequality (see \cite{D3})
\begin{Proposition}\label{Prop4.4} For any $1\le p,\rho\le+\infty$, $s\in\R$ and $0<\epsilon\le 1$,
we have
\ben\label{4.3}
\|f\|_{\widetilde{L}^\rho_T(\dot{B}^s_{p,1})}\le
C\frac{\|f\|_{\widetilde{L}^\rho_T(\dot{B}^s_{p,\infty})}}{\epsilon}
\log\bigg(e+\frac{\|f\|_{\widetilde{L}^\rho_T(\dot{B}^{s-\epsilon}_{p,\infty})}
+\|f\|_{\widetilde{L}^\rho_T(\dot{B}^{s+\epsilon}_{p,\infty})}}
{\|f\|_{\widetilde{L}^\rho_T(\dot{B}^s_{p,\infty})}}\bigg).
\een

\end{Proposition}

Now, let us prove the uniqueness of the solution of (\ref{3.1}).
Let $(u_1, h_1)$, $(u_2, h_2)$ $\in\big(L^\infty_T(B^0)\cap
L^1_T(B^2))\times L^\infty_T(\widetilde{B}^{0,1})$
be two solutions of (\ref{3.1}) with the same initial data.
The difference $\vartheta\triangleq h_2-h_1$, $w\triangleq u_2-u_1$ satisfies the following system:
\begin{align}\label{4.4}
\left\{
\begin{aligned}
&\partial_t\vartheta+u_2\cdot\nabla\vartheta=
-\dv w-w\nabla h_1-\vartheta\dv u_2-h_1\dv w,\\
&\partial_tw-\nu\widetilde{\Delta} w=-\nabla\vartheta-u_2\cdot\na w-w\cdot\nabla u_1
+\nu(1+h_1)\nabla\Big(\frac{h_1}{1+h_1}\Big)\widetilde{\nabla} w\\ &\qquad\qquad\qquad+\nu(1+h_1)
\nabla\Big(\frac{h_2}{1+h_2}-\frac{h_1}{1+h_1}\Big)\widetilde{\nabla} u_2
+\nu\vartheta\nabla\Big(\frac{h_2}{1+h_2}\Big)\widetilde{\nabla}u_2,\\
&\vartheta(0,x)=0,\quad w(0,x)=0.
\end{aligned}\right.
\end{align}
Without loss of generality, we assume that there holds for sufficiently small $T$
\begin{align}\label{4.5}
&1+h_1\ge \f12, \\\label{4.6}
&\|h_1\|_{\widetilde{E}^{0,1}_T}\le \varepsilon,
\end{align}
where $\varepsilon>0$ is small enough. Applying the Proposition \ref{Prop4.1} to the first
equation of (\ref{4.4}) yields
\begin{align}\label{4.7}
\|\vartheta(t)\|_{\dot{B}^0_{2,\infty}}\lesssim
\int_0^t\!e^{C(V_2(t)-V_2(\tau))}\|w\cdot\nabla
h_1+\vartheta\dv u_2
+h_1\dv w+\dv w\|_{\dot{B}^0_{2,\infty}}d\tau,
\end{align} with $V_2(t)\triangleq\int_0^t\|\nabla
u_2\|_{\dot{B}^1_{2,\infty}\cap L^\infty}d\tau$. It follows from
(\ref{2.5}) with $s=0$  that
\begin{align} \|w\cdot\nabla
h_1\|_{\dot{B}^0_{2,\infty}}&\lesssim
\|\nabla h_1\|_{\dot{B}^0_{2,\infty}}\|w\|_{B^1}
\lesssim\|w\|_{{B}^1}\|h_1\|_{{B}^1},\nonumber\\
\|\vartheta\dv u_2\|_{\dot{B}^0_{2,\infty}}&\lesssim\|\vartheta\|_{\dot{B}^0_{2,\infty}}
\|u_2\|_{{B}^2},\nonumber\\
\|h_1\dv w\|_{\dot{B}^0_{2,\infty}}&\lesssim
\|\dv w\|_{\dot{B}^0_{2,\infty}}
\|h_1\|_{{B}^1}\lesssim\|w\|_{{B}^1}
\|h_1\|_{{B}^1},\nonumber
\end{align}
where we have used $B^1\hookrightarrow\dot{B}^1_{2,\infty}$.
Plugging the above estimates into (\ref{4.7}), we get
\begin{align}\label{4.8}
\|\vartheta(t)\|_{\dot{B}^0_{2,\infty}}\lesssim
\int_0^t\!e^{C(V_2(t)-V_2(\tau))}\Bigl[\|w\|_{B^1}
(1+\|h_1\|_{B^1})+\|\vartheta\|_{\dot{B}^0_{2,\infty}}
\|u_2\|_{B^2}\Bigr]d\tau.
\end{align}
Recall that $u^i\in L^1_T(B^2)$, we can take a $T\in (0,\infty)$  small enough so
that
$$C\|u_2\|_{L^1_T(B^2)}\le \f1 4,$$
which together with (\ref{4.8}) implies that for $t\le T$
\begin{align}\label{4.9}
\|\vartheta\|_{L^\infty_t(\dot{B}^0_{2,\infty})}\lesssim
\|w\|_{L^1_t(B^1)}(1+\|h_1\|_{L^\infty_t(B^1)}).
\end{align}
Applying (\ref{4.3}) to the term
$\|w\|_{L^1_t(B^1)}$ yields
\begin{align}\label{4.10}
\|\vartheta\|_{L^\infty_t(\dot{B}^0_{2,\infty})}\lesssim&
\|w\|_{\widetilde{L}^1_t(\dot{B}^1_{2,\infty})}
\log\bigg(e+\frac{\|w\|_{\widetilde{L}^1_t(\dot{B}^0_{2,\infty})}
+\|w\|_{\widetilde{L}^1_t(\dot{B}^2_{2,\infty})}}
{\|w\|_{\widetilde{L}^1_t(\dot{B}^1_{2,\infty})}}\bigg)(1+\|h_1\|_{L^\infty_t(B^1)}).
\end{align}
Thanks to $B^s\hookrightarrow\dot{B}^s_{2,\infty}$
and $\widetilde{B}^{0,1}\hookrightarrow B^1$, we have
\begin{align}\label{4.11}
\|\vartheta\|_{L^\infty_t(\dot{B}^0_{2,\infty})}\lesssim
\|w\|_{\widetilde{L}^1_t(\dot{B}^1_{2,\infty})}
\log\bigg(e+\frac{W(t)}
{\|w\|_{\widetilde{L}^1_t(\dot{B}^1_{2,\infty})}}\bigg)
\end{align}
with $$W(t)\triangleq\|u^i\|_{\widetilde{L}^1_t(B^0)}
+\|u^i\|_{\widetilde{L}^1_t(B^2)},$$ and for finite $t$, $W(t)< +\infty$.

Next, we deal with the second equation of (\ref{4.4}). We get by
applying (\ref{2.7}) with $s=1$, $s=0$ respectively that
\begin{align}\label{4.12}
\|u_2\cdot\nabla w\|_{\widetilde{L}^1_t(\dot{B}^{-1}_{2,\infty})}\lesssim
\|u_2\|_{L^2_t(B^1)}\|w\|_{\widetilde{L}^2_t(\dot{B}^{0}_{2,\infty})},\\\label{5.12}
\|w\cdot\nabla u_1\|_{\widetilde{L}^1_t(\dot{B}^{-1}_{2,\infty})}\lesssim
\|u_1\|_{L^2_t(B^1)}\|w\|_{\widetilde{L}^2_t(\dot{B}^{0}_{2,\infty})}.
\end{align}
We can deduce $h_i\in C(0,T; \R^2)$ $(i=1,2)$ from the fact $B^1\hookrightarrow C$.
Moreover, due to (\ref{4.5}), we can assume $h_1(t,x)+1\ge \delta$
for all $t\le T$, $x\in \R^2$. Since $h_1$, $h_2$ have the same initial data,
from the continuity of $h_2$, there exists a $\widetilde{T}\le T$ such that
$$h_2(x,t)+1\ge \delta,\quad\hbox{for all}\quad t\in[0, \widetilde{T}],\quad x\in\R^2.$$
It follows from (\ref{2.6}) with $s=1$, (\ref{2.13}) and
$B^1\hookrightarrow\dot{B}^1_{2,\infty}\cap L^\infty$ that
\begin{align}
&\Big\|(1+h_1)\nabla\Big(\frac{h_2}{1+h_2}-\frac{h_1}{1+h_1}\Big)\widetilde{\nabla} u_2\Big\|
_{\dot{B}^{-1}_{2,\infty}}\nonumber\\&\qquad\lesssim
\Big\|(1+h_1)\nabla\Big(\frac{h_2}{1+h_2}-\frac{h_1}{1+h_1}\Big)\Big\|
_{\dot{B}^{-1}_{2,\infty}}\|\widetilde{\nabla} u_2\|
_{B^1}\nonumber\\&\qquad\lesssim (1+\|h_1\|_{B^1})
\Big\|\frac{h_2}{1+h_2}-\frac{h_1}{1+h_1}\Big\|
_{\dot{B}^0_{2,\infty}}\|u_2\|_{B^2}\nonumber\\&\qquad\lesssim (1+\|h_1\|_{B^1})
(\|h_1\|_{B^1}+\|h_2\|_{B^1})\|\vartheta\|_{\dot{B}^0_{2,\infty}}
\|u_2\|_{B^2},\nonumber
\end{align}
which together with
 $L^1_t(\dot{B}^{-1}_{2,\infty})\subset\widetilde{L}^1_t(\dot{B}^{-1}_{2,\infty})$
yields\begin{align}\label{4.14}
\Big\|&(1+h_1)\nabla\Big(\frac{h_2}{1+h_2}-\frac{h_1}{1+h_1}\Big)\widetilde{\nabla} u_2\Big\|
_{\widetilde{L}^1_t(\dot{B}^{-1}_{2,\infty})}\nonumber\\&
\lesssim\int_0^t(1+\|h_1\|_{B^1})(\|h_1\|_{B^1}+\|h_2\|_{B^1})
\|\vartheta\|_{\dot{B}^0_{2,\infty}}\|u_2\|_{B^2}d\tau.
\end{align}
Thanks to (\ref{2.6}),
(\ref{2.12}),   and
$L^1_t(\dot{B}^{-1}_{2,\infty})\subset\widetilde{L}^1_t(\dot{B}^{-1}_{2,\infty})$, we get
\begin{align}\label{4.15}
&\Big\|\vartheta\nabla \Big(\frac{h_2}{1+h_2}\Big)\widetilde{\nabla} u_2\Big\|
_{\widetilde{L}^1_t(\dot{B}^{-1}_{2,\infty})}
\lesssim\int_0^t\|\vartheta\|_{\dot{B}^{0}_{2,\infty}}\|h_2\|_{B^1}\|u_2\|_{B^2}d\tau.\end{align}
Thanks to Lemma 5.3 with $s_1=s_2=0$,  Lemma 5.4 with $s=1$, and (\ref{2.5})  with $s=0$, we have
\begin{align} \label{4.17}
\sup_{k\in\Z}&\, \omega_k(t)2^{-k}\Big\|\Delta_k\Big((1+h_1)\nabla\Big(\frac{h_1}{1+h_1}\Big)\widetilde{\nabla} w\Big)\Big\|_{L^1_t(L^2)}
\nonumber\\&\lesssim
\Big\|\na\Big(\frac{h_1}{1+h_1}\Big)\Big\|_{E^{0}_t}
\|(1+h_1)\widetilde{\nabla} w\|_{\widetilde{L}^1_t(\dot{B}^0_{2,\infty})}\nonumber\\&\lesssim
\Big\|\frac{h_1}{1+h_1}\Big\|_{E^{1}_t}(1+\|h_1\|_{L^\infty_t(B^1)})
\|\nabla w\|_{\widetilde{L}^1_t(\dot{B}^0_{2,\infty})}
\nonumber\\&\lesssim\|h_1\|_{\widetilde{E}^{0,1}_t}
(1+\|h_1\|_{L^\infty_t(\widetilde{B}^{0,1})})^4\|w\|_{\widetilde{L}^1_t(\dot{B}^1_{2,\infty})}.
\end{align}
In terms of  Proposition 4.3, (\ref{4.12})-(\ref{4.17}) and
$\widetilde{B}^{0,1}\hookrightarrow B^1$, we finally obtain
\begin{align}\label{4.18}
\|w\|&_{\widetilde{L}^1_t(\dot{B}^{1}_{2,\infty})}+
\|w\|_{\widetilde{L}^2_t(\dot{B}^{0}_{2,\infty})}\nonumber\\ \lesssim &
\|u_2\|_{L^2_t(B^1)}\|w\|_{\widetilde{L}^2_t(\dot{B}^{0}_{2,\infty})}
+\|u_1\|_{L^2_t(B^1)}\|w\|_{\widetilde{L}^2_t(\dot{B}^{0}_{2,\infty})}\nonumber\\&+\|h_1\|_{\widetilde{E}^{0,1}_t}
\big(1+\|h_1\|_{L^\infty_t(\widetilde{B}^{0,1})}\big)^4
\|w\|_{\widetilde{L}^1_t(B^1_{2,\infty})}\nonumber\\
&+\int_0^t(1+\|h_1\|_{\widetilde{B}^{0,1}})(1+\|h_1\|_{\widetilde{B}^{0,1}}
+\|h_2\|_{\widetilde{B}^{0,1}})(1+\|u_2\|_{B^2})
\|\vartheta\|_{\dot{B}^0_{2,\infty}}d\tau.
\end{align}
Let us define $$Z(t)\triangleq \|w\|_{\widetilde{L}^1_t(\dot{B}^{1}_{2,\infty})}
+\|w\|_{\widetilde{L}^2_t(\dot{B}^{0}_{2,\infty})}.$$
Due to (\ref{4.6}), if $T$ is chosen small enough,
then the first four terms of the right side of (\ref{4.18}) can be absorbed by the left
side $Z(t)$. Noting that $r\log(e+\frac{W(T)}{r})$ is increasing, from
(\ref{4.11}) and (\ref{4.18}), it follows that
\begin{align}\label{4.19}
Z(t)&\lesssim\int_0^t(1+W'(\tau))
Z(\tau)\log\Big(e+\frac{W(\tau)}{Z(\tau)}\Big)d\tau\nonumber\\&
\lesssim \int_0^t(1+W'(\tau))
Z(\tau)\log\Big(e+\frac{W(T)}{Z(\tau)}\Big)d\tau.
\end{align}
It is easy to verify that
$$1+W'(\tau)\in L^1_{loc}(\R^+)\quad\hbox{and}
\quad\int_0^1\frac{dr}{r\log(e+\frac{W(T)}{r})}=+\infty.$$
Hence by Osgood Lemma, we have $Z\equiv0$ on $[0,\widetilde{T}]$, i.e. $w\equiv0$, then from
(\ref{4.9}),
$\vartheta=h_2-h_1\equiv0$. Then a standard
continuous argument gives the uniqueness.

\setcounter{equation}{0}
\section{Appendix}

In this appendix, we prove some  multilinear estimates in the weighted Besov space.

\begin{Lemma}\label{Lem5.1}
Let $A$ be a homogeneous smooth function of degree $m$. Assume that  $-\frac d2<\rho\le\frac d2$. Then there hold
\ben
&&\Big|(A(D)\Delta_k(v\cdot\na h),A(D)\Delta_kh)\Big|\le C\|{\cal
F}_k^m(t)\|_2\|A(D)\Delta_kh\|_2,\label{5.1a}\\
&&\Big|(A(D)\Delta_k(v\cdot\na u),A(D)\Delta_ku)\Big|\le C\|{\cal
\widetilde{F}}_k^m(t)\|_2\|A(D)\Delta_ku\|_2,\label{5.2a}
\een
and
\ben
&&\Big|(A(D)\Delta_k(v\cdot\na h),\Delta_k u)+(\Delta_k(v\cdot\na
u),A(D)\Delta_kh)\Big|\nonumber\\
&&\qquad\qquad\le C\big(\|{\cal {F}}_k^m(t)\|_2+\|{\cal \widetilde{F}}_k^0(t)\|_2\big)\big(\|\Delta_k
u\|_2+\|A(D) \Delta_kh\|_2\big),\label{5.3a}
\een
where ${\cal F}_k^m(t)$ and ${\cal \widetilde{F}}_k^m(t)$ satisfy
\ben\label{5.1}
&&\sum_{k\in\Z}\omega_k(T)2^{k(\rho-m)}\|{\cal
F}_k^m(t)\|_{L^1_T(L^2)}\le C\|h\|_{E^\rho_T}\|v\|_{L^1_T(B^{\frac{d}2+1})},\\
&&\sum_{k\in\Z}2^{k(\rho-m)}\|{\cal \widetilde{F}}_k^m(t)\|_{L^1_T(L^2)}
\le C\|u\|_{L^2_T(B^{\rho+1})}\|v\|_{L^2_T(B^{\frac{d}2})}.\label{5.2}
\een
\end{Lemma}

\no{\it Proof.}\quad Let us first prove (\ref{5.1a}).
Using the Bony's paraproduct decomposition, we write
\begin{align}\label{5.3}
\big(A(D)\Delta_k(v\cdot\na h), A(D)\Delta_kh\big)&=\big(A(D)\Delta_k(T'_{\pa_j h}v^j),  A(D)\Delta_kh\big)+J_k,
\end{align}
where
\beno
&&T'_fg=T_fg+R(f,g),\quad \textrm{and} \\
&&J_k=\sum_{|k'-k|\le3}([A(D)\Delta_k, S_{k'-1}v^j]\Delta_{k'}\pa_jh,  A(D)\Delta_kh)\\
&&\qquad+\sum_{|k'-k|\le3}((S_{k'-1}-S_{k-1})v^jA(D)\Delta_{k}\Delta_{k'}\pa_jh,  A(D)\Delta_kh)\\
&&\qquad+(S_{k-1}v^jA(D)\Delta_{k}\pa_jh,  A(D)\Delta_kh)
\eeno
We get by integration by parts that
$$(S_{k-1}v^jA(D)\Delta_{k}\pa_jh,  A(D)\Delta_kh)=-\frac12\big(S_{k-1}\dv\, vA(D)\Delta_{k}h,  A(D)\Delta_kh\big).$$
Let us set
\beno
&&{\cal F}_{k,0}^m(t)=A(D)\Delta_k(T'_{\pa_j h}v^j),\\
&&{\cal F}_{k,1}^m(t)=\sum_{|k'-k|\le3}[A(D)\Delta_k, S_{k'-1}v^j]\Delta_{k'}\pa_jh,\\
&&{\cal F}_{k,2}^m(t)=\sum_{|k'-k|\le3}(S_{k'-1}-S_{k-1})v^jA(D)\Delta_{k}\Delta_{k'}\pa_jh,\\
&&{\cal F}_{k,4}^m(t)=-\frac12 S_{k-1}\dv\, vA(D)\Delta_{k}h.
\eeno
By the Cauchy-Schwartz inequality, we get
\beno
\big(A(D)\Delta_k(v\cdot\na h), A(D)\Delta_kh\big)\le \|{\cal
F}_k^m(t)\|_2\|A(D)\Delta_kh\|_2,
\eeno
with ${\cal F}_k^m(t)=\ds\sum_{i=0}^3{\cal F}_{k,i}^m(t).$ So,
it remains to prove that ${\cal F}_k^m(t)$ satisfies (\ref{5.1}). For the simplicity, we set
$$\widetilde{\Delta}_{k}=\ds\sum_{|k'-k|\le 1}\Delta_{k'},\quad
\widetilde{\widetilde{\Delta}}_{k}=\ds\sum_{|k'-k|\le 3}\Delta_{k'}.$$
Thanks to the definition of ${\cal F}_{k,0}^m(t)$ and Lemma \ref{Lem2.1}, we have
\begin{align}
\|{\cal F}_{k,0}^m(t)\|_{L^1_T(L^2)}&\le \sum_{|k'-k|\le3}2^{km}\|S_{k'-1}\pa_jh\|_{L^\infty_T(L^\infty)}\|\Delta_{k'}v^j\|_{L^1_T(L^2)}\nonumber\\
&\quad+\sum_{k'\ge k-2}2^{(m+\frac d2)k}\|\Delta_k(\Delta_{k'}\pa_jh\widetilde{\Delta}_{k'}v^j)\|_{L^1_T(L^1)}\nonumber\\
&\triangleq I+II.\nonumber
\end{align}
Thanks to Lemma \ref{Lem2.1}, we have
\begin{align}
2^{k(\rho-m)}I &\lesssim 2^{k\rho}
\sum_{k'\le k+1}2^{k'(1+\f d2)}\|\Delta_{k'}h\|_{L^\infty_T(L^2)}\|\widetilde{\widetilde{\Delta}}_{k}v\|_{L^1_T(L^2)}\nonumber\\
&\lesssim\sum_{k'\le k+1}2^{(k'-k)(1+\f d2-\rho)}2^{k'\rho}\|\Delta_{k'}h\|_{L^\infty_T(L^2)}2^{k(1+\f d2)}\|\widetilde{\widetilde{\Delta}}_{k}v\|_{L^1_T(L^2)},\nonumber
\end{align}
from which and the definition of $\omega_k(T)$, it follows that
\ben\label{5.4}
&&\sum_{k\in \Z}\omega_k(T)2^{k(\rho-m)}I\nonumber\\
&&\lesssim \sum_{k'\in \Z}2^{k'\rho}\|\Delta_{k'}h\|_{L^\infty_T(L^2)}\sum_{k\ge k'-1}\omega_k(T)2^{(k'-k)(1+\f d2-\rho)}2^{k(1+\f d2)}\|\widetilde{\widetilde{\Delta}}_{k}v\|_{L^1_T(L^2)}\nonumber\\
&&\lesssim \sum_{k'\in \Z}\omega_{k'}(T)2^{k'\rho}\|\Delta_{k'}h\|_{L^\infty_T(L^2)}\sum_{k\ge k'-1}2^{(k'-k)(\f d2-\rho)}2^{k(1+\f d2)}\|\widetilde{\widetilde{\Delta}}_{k}v\|_{L^1_T(L^2)}\nonumber\\
&&\lesssim \|h\|_{E^\rho_T}\|v\|_{L^1_T(B^{\f d2+1})}, \een where we
used the assumption $\rho\le \f d 2$ in the last inequality. Set
$e_k(T)=e^1_{k}(T)+e^2_{k}(T)$. Using Lemma \ref{Lem2.1}, we also have
\begin{align}
\omega_k(T)2^{k(\rho-m)}II&\lesssim \omega_k(T)2^{k(\rho+\f d2)}
\sum_{k'\ge k-2}2^{k'}\|\Delta_{k'}h\|_{L^\infty_T(L^2)}\|\widetilde{\Delta}_{k'}v\|_{L^1_T(L^2)}\nonumber\\
&\lesssim 2^{k(\rho+\f d2)} \sum_{k'\ge
k-2}2^{k'}\|\Delta_{k'}h\|_{L^\infty_T(L^2)}\|\widetilde{\Delta}_{k'}v\|_{L^1_T(L^2)}
\sum_{k'\ge \widetilde{k}\ge k}2^{-(\widetilde{k}-k)}e_{\widetilde{k}}(T)\nonumber\\
&\quad+2^{k(\rho+\f d2)} \sum_{k'\ge
k-2}2^{k'}\|\Delta_{k'}h\|_{L^\infty_T(L^2)}\|\widetilde{\Delta}_{k'}v\|_{L^1_T(L^2)}
\sum_{\widetilde{k}\ge k, \widetilde{k}\ge k'}2^{-(\widetilde{k}-k)}e_{\widetilde{k}}(T)\nonumber\\
&\triangleq II_1+II_2.\nonumber
\end{align}
Note that for $\widetilde{k}\le k'$
$$
e_{\widetilde{k}}(T)\le e_{k'}(T)\le \omega_{k'}(T),
$$
from which and $\rho>-\f d 2$, we deduce that
\begin{align}\label{5.5}
\sum_{k\in Z}II_1&\lesssim \sum_{k'\in Z}\omega_{k'}(T)2^{k'\rho}\|\Delta_{k'}h\|_{L^\infty_T(L^2)}2^{k'(\f d2+1)}
\|\widetilde{\Delta}_{k'}v\|_{L^1_T(L^2)}
\sum_{k\le k'+2}2^{(k-k')(\rho+\f d2)}\nonumber\\
&\lesssim \|h\|_{E^\rho_T}\|v\|_{L^1_T(B^{\f d2+1})}.
\end{align}
Similarly, we can obtain
\begin{align}\label{5.6}
\sum_{k\in \Z}II_2&\lesssim \sum_{k'\in \Z}2^{k'\rho}\|\Delta_{k'}h\|_{L^\infty_T(L^2)}\sum_{k\le k'+2}2^{(k-k')(\f d2+\rho)}
\sum_{\widetilde{k}\ge k'}2^{-(\widetilde{k}-k)}e_{\widetilde{k}}(T)\|v\|_{L^1_T(B^{\f d2+1})}\nonumber\\
&\lesssim \sum_{k'\in \Z}\omega_{k'}(T)2^{k'\rho}\|\Delta_{k'}h\|_{L^\infty_T(L^2)}\sum_{k\le k'+2}2^{(k-k')(\f d2+\rho+1)}
\|v\|_{L^1_T(B^{\f d2+1})}\nonumber\\
&\lesssim  \|h\|_{E^\rho_T}\|v\|_{L^1_T(B^{\f d2+1})}.
\end{align}

By summing up (\ref{5.4})-(\ref{5.6}), we obtain
\ben\label{5.7}
\sum_{k\in\Z}\omega_k(T)2^{k(\rho-m)}\|{\cal F}_{k,0}^m(t)\|_{L^1_T(L^2)}\lesssim \|h\|_{E^\rho_T}\|v\|_{L^1_T(B^{\f d2+1})}.
\een

Note that $A(D)\Delta_k=2^{km}\widetilde{\varphi}(2^{-k}D)$ with $\widetilde{\varphi}(\xi)=A(\xi)\varphi(\xi)$.
Set $\tilde\theta={\cal F}^{-1}\tilde\varphi$, we get by using the Taylor's formula that
\begin{align}\label{5.8}
{\cal F}_{k,1}^m(t)=\sum_{|k'-k|\le3}2^{k(m-1)}\int_{\R^d}\int_0^1\tilde\theta(y)
(y\cdot S_{k'-1}\na v^j(x-2^{-k}\tau y))\Delta_{k'}\pa_jh(x-2^{-k}y) d\tau dy,\nonumber
\end{align}
from which and Lemma \ref{Lem2.1}, it follows that
\begin{align}
\|{\cal F}_{k,1}^m(t)\|_{L^1_T(L^2)}&\lesssim2^{k(m-1)}\sum_{|k'-k|\le3}
\|S_{k'-1}\na v^j\|_{L^1_T(L^\infty)}\|\Delta_{k'}\pa_jh\|_{L^\infty_T(L^2)}\nonumber\\&
\lesssim2^{km}\sum_{|k'-k|\le3}
\|\Delta_{k'}h\|_{L^\infty_T(L^2)}\|v\|_{L^1_T(B^{\frac d 2 +1})},\nonumber
\end{align}
thus, we get
\begin{align}
\sum_{k\in\Z}\omega_k(T)2^{k(\rho-m)}\|{\cal F}_{k,1}^m(t)\|_{L^1_T(L^2)}
\lesssim  \|h\|_{E^\rho_T}\|v\|_{L^1_T(B^{\f d2+1})}.
\end{align}
Thanks to  the fact $|k'-k|\le3$ and Lemma \ref{Lem2.1}, we have
\begin{align}
\|(S_{k'-1}-S_{k-1})v^jA(D)\Delta_{k}\Delta_{k'}\pa_jh\|_{L^1_T(L^2)}
\lesssim 2^{km}\|\Delta_{k}h\|_{L^\infty_T(L^2)}\|v\|_{L^1_T(B^{\f d2+1})},\nonumber
\end{align}
from which, it follows that
\begin{align}\label{5.9}
\sum_{k\in\Z}&\omega_k(T)2^{k(\rho-m)}\big(\|{\cal F}_{k,2}^m(t)\|_{L^1_T(L^2)}
+\|{\cal F}_{k,3}^m(t)\|_{L^1_T(L^2)}\big)\lesssim\|h\|_{E^\rho_T}\|v\|_{L^1_T(B^{\f d2+1})},
\end{align}
which together with (\ref{5.7}) and (\ref{5.8}) yields (\ref{5.1}).

Using the decomposition (\ref{5.3}) with $h$ instead of $u$ and Lemma \ref{Lem2.1}, (\ref{5.2a}) can be easily proved. We omit it here.
In order to prove (\ref{5.3a}), we use the decomposition
\begin{align}
\big(A(D)\Delta_k(v\cdot\na h),\Delta_k u\big)+\big(\Delta_k(v\cdot\na u),A(D)\Delta_kh\big)=I_k+J_k,\nonumber
\end{align}
with
\begin{align}
I_k=&\big(A(D)\Delta_k(T'_{\pa_j h}v^j),\Delta_k u\big)+\big(\Delta_k(T'_{\pa_j u}v^j),A(D)\Delta_kh\big)\nonumber\\
\triangleq&\big({\cal F}_{k,0}^m(t),\Delta_k u\big)+\big({\cal \widetilde{F}}_{k,0}^0(t),A(D)\Delta_kh\big)\nonumber\\
J_k=&\sum_{|k'-k|\le3}\Big([A(D)\Delta_k, S_{k'-1}v^j]\Delta_{k'}\pa_jh,  \Delta_ku\Big)+
\Big((S_{k'-1}-S_{k-1})v^jA(D)\Delta_{k}\Delta_{k'}\pa_jh,  \Delta_ku\Big)\nonumber\\
&+\sum_{|k'-k|\le3}\Big([\Delta_k, S_{k'-1}v^j]\Delta_{k'}\pa_ju,  A(D)\Delta_kh\Big)+
\Big((S_{k'-1}-S_{k-1})v^j\Delta_{k}\Delta_{k'}\pa_ju,  A(D)\Delta_kh\Big)\nonumber\\
&-\Big(S_{k-1}\dv vA(D)\Delta_{k}h,  \Delta_ku\Big)\nonumber\\
\triangleq&\big({\cal F}_{k,1}^m(t),\Delta_k u\big)+\big({\cal F}_{k,2}^m(t),\Delta_k u\big)+
\big({\cal \widetilde{F}}_{k,1}^0(t),A(D)\Delta_kh\big)
+\big({\cal\widetilde{ F}}_{k,2}^0,A(D)\Delta_kh\big)+\big({\cal F}_{k,3}^m(t), \Delta_ku\big),
\nonumber\end{align}
from which, a similar proof of (\ref{5.1}) gives (\ref{5.3a}). This completes the proof of Lemma \ref{Lem5.1}.\ef

\begin{Lemma}\label{Lem5.2}
Let  $s_1\le\frac d2-1$, $s_2\le\frac d2$, and $s_1+s_2>0$. Then there holds
\begin{align}
\sum_{k\in\Z}\omega_k(T)2^{k(s_1+s_2-\frac d2)}\|\Delta_k(fg)\|_{L^1_T(L^2)}
\le C\sum_{k\in\Z}\omega_k(T)2^{ks_1}\|\Delta_kf\|_{L^{r_1}_T(L^2)}\|g\|_{L^{r_2}_T(B^{s_2})},\label{5.10}
\end{align}
where $1\le r_1,r_2\le \infty$ and $\f 1 {r_1}+\f 1 {r_2}=1$.
\end{Lemma}

\no{\it Proof.}\quad Using the Bony's paraproduct decomposition, we write
\begin{align}
\Delta_k(fg)&=\sum_{|k'-k|\le3}\Delta_k(S_{k'-1}f\Delta_{k'}g)+
\sum_{|k'-k|\le3}\Delta_k(S_{k'-1}g\Delta_{k'}f)\nonumber\\
&\quad+\sum_{k'\ge k-2}\Delta_k(\Delta_{k'}f\widetilde{\Delta}_{k'}g)\triangleq I+II+III.\nonumber
\end{align}
A similar proof of (\ref{5.4}) ensures that for $s_1\le\frac d2-1$
\begin{align}
\sum_{k\in\Z}\omega_k(T)2^{k(s_1+s_2-\frac d2)}\|I\|_{L^1_T(L^2)}\lesssim\sum_{k\in\Z}\omega_k(T)2^{ks_1}\|\Delta_kf\|_{L^{r_1}_T(L^2)}\|g\|_{L^{r_2}_T(B^{s_2})},\nonumber
\end{align}
while $II$ can be directly deduced for $s_2\le\frac d2$. On the other hand,
a similar proof of (\ref{5.5}) and (\ref{5.6}) gives for $s_1+s_2>0$
\begin{align}
\sum_{k\in\Z}\omega_k(T)2^{k(s_1+s_2-\frac d2)}\|III\|_{L^1_T(L^2)}\lesssim \sum_{k\in\Z}\omega_k(T)2^{ks_1}\|\Delta_kf\|_{L^{r_1}_T(L^2)}\|g\|_{L^{r_2}_T(B^{s_2})}.\nonumber
\end{align}
This completes the proof of Lemma \ref{Lem5.2}.\ef

Similarly, we can also prove the following lemma.

\begin{Lemma}\label{Lem5.3}
Let $s_1\le\frac d2-1$, $s_2<\frac d2$ and $s_1+s_2\ge0$. Then there holds
\begin{align}
\sup_{k\in\Z}\omega_k(T)2^{k(s_1+s_2-\frac d2)}\|\Delta_k(fg)\|_{L^1_T(L^2)}
\le C\|f\|_{E^{s_1}_T}\|g\|_{\widetilde{L}^1_T(\dot B^{s_2}_{2,\infty})}.\label{5.11}
\end{align}
\end{Lemma}

\begin{Lemma}\label{Lem5.4}
Let $s>0$. Assume that $F\in W^{[s]+3,\infty}_{loc}(\R^d)$ with  $F(0)=0$. Then
there holds
\ben
\|F(f)\|_{E^s_T}\le C(1+\|f\|_{L^\infty_T(L^\infty)})^{[s]+2}\|f\|_{E^s_T}.\label{5.12}
\een
\end{Lemma}

\no{\it Proof.}\quad We decompose $F(f)$ as
\begin{align}
F(f)=\sum_{k'\in\Z}F(S_{k'+1}f)-F(S_{k'}f)&=\sum_{k'\in\Z}\Delta_{k'}f\int_0^1F'(S_{k'}f+\tau\Delta_{k'}f)d\tau\nonumber\\
&\triangleq\sum_{k'\in\Z}\Delta_{k'}f\, m_{k'},\nonumber
\end{align}
where $m_{k'}=\int_0^1F'(S_{k'}f+\tau\Delta_{k'}f)d\tau$. Furthermore, we write
\begin{align}
\Delta_kF(f)=\sum_{k'<k}\Delta_k(\Delta_{k'}f\,m_{k'})+\sum_{k'\ge k}\Delta_k(\Delta_{k'}f\, m_{k'})\triangleq I+II.\nonumber
\end{align}
By Lemma \ref{Lem2.1}, we have
\begin{align}\label{5.13}
\|I\|_{L^\infty_T(L^2)}&\le\sum_{k'<k}\|\Delta_k(\Delta_{k'}f\,m_{k'})\|_{L^\infty_T(L^2)}\nonumber\\
&\le \sum_{k'<k}2^{-k|\al|}\sup_{|\gamma|=|\al|}\|D^\gamma\Delta_k(\Delta_{k'}f\,m_{k'})\|_{L^\infty_T(L^2)},
\end{align}
with $\al$ to be determined later. Note that for $|\gamma|\ge0$, we have
$$\|D^\gamma m_{k'}\|_\infty\lesssim2^{k'|\gamma|}(1+\|f\|_\infty)^{|\gamma|}\|F'\|_{W^{|\gamma|,\infty}},$$
from which and (\ref{5.13}), it follows that
\begin{align}
2^{ks}\|I\|_{L^\infty_T(L^2)}\lesssim2^{k(s-|\al|)}\sum_{k'<k}2^{k'|\al|}
\|\Delta_{k'}f\|_{L^\infty_T(L^2)}(1+\|f\|_{L^\infty_T(L^\infty)})^{|\al|}\|F'\|_{W^{|\al|,\infty}},\nonumber
\end{align}
thus, if we take $|\al|=[s]+2$, we get
\begin{align}
&\sum_{k\in\Z}\omega_k(T)2^{ks}\|I\|_{L^\infty_T(L^2)}\nonumber\\
&\quad\lesssim\sum_{k'\in\Z}2^{k's}\omega_{k'}(T)\|\Delta_{k'}f\|_{L^\infty_T(L^2)}
\sum_{k>k'}2^{(k-k')(s-|\al|+1)}(1+\|f\|_{L^\infty_T(L^\infty)})^{|\al|}\|F'\|_{W^{|\al|,\infty}}
\nonumber\\&\quad \lesssim(1+\|f\|_{L^\infty_T(L^\infty)})^{[s]+2}\|F'\|_{W^{[s]+2,\infty}}\|f\|_{E^s_T}.\label{5.14}
\end{align}

Next, let us turn to the proof of $II$. We get by using Lemma 2.1 that
\begin{align}
\|II\|_{L^\infty_T(L^2)}&\lesssim\sum_{k\ge k'}\|\Delta_{k'}f\|_{L^\infty_T(L^2)}.\nonumber
\end{align}
Then we write
\begin{align}
\sum_{k\in\Z}\omega_k(T)2^{ks}\|II\|_{L^\infty_T(L^2)}
&\lesssim \sum_{k\in \Z}2^{ks}\sum_{k'\ge k}\|\Delta_{k'}f\|_{L^\infty_T(L^2)}
\sum_{k'\ge \widetilde{k}\ge k}2^{-(\widetilde{k}-k)}e_{\widetilde{k}}(T)\nonumber\\
&\quad+\sum_{k\in \Z}2^{ks}\sum_{k'\ge k}\|\Delta_{k'}f\|_{L^\infty_T(L^2)}
\sum_{\widetilde{k}\ge k', \widetilde{k}\ge k}2^{-(\widetilde{k}-k)}e_{\widetilde{k}}(T),\nonumber
\end{align}
from which, a similar proof of (\ref{5.5}) and (\ref{5.6}) ensures that
\begin{align}
\sum_{k\in\Z}\omega_k(T)2^{ks}\|II\|_{L^\infty_T(L^2)}\lesssim\|f\|_{E^s_T}.\label{5.15}
\end{align}

By summing up (\ref{5.14}) and (\ref{5.15}), we deduce the inequality (\ref{5.12}).
This completes the proof of Lemma \ref{Lem5.4}.\ef

\bigskip

\noindent {\bf Acknowledgments.}  The authors thank the referees
 for their invaluable comments and suggestions which helped improve the paper greatly.
  Q. Chen and C.
Miao  were partially  supported by the NSF of China (No.10571016).
Z. Zhang is supported by the NSF of China(No.10601002).


\begin{thebibliography}{50}
\bibitem{BC} H. Bahouri, J.-Y. Chemin,
{\it Equations de transport relatives \`{a} des champs des
vecteurs non-lipschitziens et m\'{e}canique des fluides}. Arch.
Rational Mech. Anal., 127(1994), 159-199.

\bibitem{BC1} H. Bahouri, J.-Y. Chemin,
{\it \'{E}quations d'ondes quasilin\'{e}aires et estimations de Strichartz}.
Amer. J. Math., 121(1999), 1337-1377.


\bibitem{B} J.-M. Bony, {\it Calcul symbolique et propagation
des singulariti\'{e}s pour les \'{e}quations aux d\'{e}riv\'{e}es
partielles non lin\'{e}aires}. Ann. de l'Ecole Norm. Sup.,
14(1981), 209-246.

\bibitem{BDM} D. Bresch, B. Desjardins, G. M\'{e}tivier, {\it Recent mathematical results and open problems about shallow water equations},
Preprint.

\bibitem{Bui} A. T. Bui, {\it Existence and uniqueness of a classical solution of an
initial boundary value problem of the theroy of shallow waters}.
SIAM J. Math. Anal., 12(1981), 229-241.


\bibitem{Can1}M. Cannone, {\it Ondelettes, paraproduits et Navier-Stokes}.
Nouveaux essais, Diderot \'{e}diteurs, Paris, 1995.


\bibitem{Can2}M. Cannone, {\it Harmonic analysis tools for solving the incompressible
                     Navier-Stokes equations}.
Handbook of Mathematical fluid Dynamics, Vol. III, North-Holland, Amsterdam, 2004.

\bibitem{Ch1} J.-Y. Chemin,
{\it Perfect incompressible fluids}. Oxford University Press, New York, 1998.


\bibitem{Ch2} J.-Y. Chemin,
{\it Th\'{e}or\`{e}mes d'unicit\'{e} pour le syst\`{e}me de Navier-Stokes tridimensionnel}.
J. d'Analyse Math., 77(1999), 27-50.



\bibitem{CL} J.-Y. Chemin, N. Lerner,
{\it Flot de champs de vecteurs non lipschitziens et \'{e}quations de
Navier-Stokes}. J. Differential Equations, 121(1992), 314-328.

\bibitem{D1} R. Danchin,
{\it Global existence  in  critical spaces for compressible Navier-Stokes equations}.
Invent. Math., 141(2000), 579-614.

\bibitem{D2} R. Danchin,
{\it Global existence in  critical spaces for flows of compressible
viscous and heat-conductive gases}.
Arch.
Rational Mech. Anal., 160(2001), 1-39.

\bibitem{D5} R. Danchin,
{\it Local theory in critical spaces for compressible viscous and heat-conductive gases},
Comm. Partial Differential Equations, 26(2001),1183-1233.

\bibitem{D3} R. Danchin,
{\it Density-dependent incompressible viscous fluids in critical spaces}.
Proc. Roy. Soc. Edinburgh Sect. A, 133(2003), 1311-1334.


\bibitem{D4} R. Danchin,
{\it On the uniqueness  in critical spaces for compressible Navier-Stokes equations}.
Nonlinear Differential Equations Appl., 12(2005), 111-128.



\bibitem{FK} H. Fujita, T. Kato,
{\it On the Navier-Stokes initial value problem I}.
Arch.
Rational Mech. Anal., 16(1964), 269-315.

\bibitem{Klo} P. E. Kloeden,
{\it Global existence of classical solutions in the dissipative shallow water equations}.
SIAM J. Math. Anal., 16(1985), 301-315.

\bibitem{Mey2} Y. Meyer, {\it Wavelets, paraproducts  and Navier-Stokes equations.
Current developments in mathematics}. International Press, 1996.


\bibitem{RS} T. Runst, W. Sickel,
{\it Sobolev spaces of fractional order, Nemytskij operators, and
nonlinear partial differential Equations}. de Gruyter Series in
Nonlinear Analysis and Applications, 3. Walter de Gruyter \& Co.,
Berlin, 1996.

\bibitem{Su1} L. Sundbye,
{\it Global existence for the Dirichlet problem for the viscous shallow water equations}.
J. Math. Anal. Appl., 202(1996), 236-258.

\bibitem{Su2} L. Sundbye,
{\it Global existence for the Cauchy problem for the viscous shallow water equations}.
Rocky Mountain J. Math., 28(1998), 1135-1152.


\bibitem{Tr} H. Triebel,
{\it Theory of Function Spaces}. Monographs in Mathematics,
 Birkh\"{a}user Verlag, Basel, Boston, Stuttgart, 1983.

\bibitem{WX} W-K. Wang, C-J. Xu,
{\it The Cauchy problem for viscous shallow water equations}.
Rev. Mat. Iberoamericana, 21(2005), 1-24.





\end{thebibliography}
\end{document}